\documentclass[11pt]{article}
\usepackage[english]{babel}
\usepackage{times}
\usepackage[T1]{fontenc}
\usepackage{amsmath, amsthm, amsfonts}
\usepackage{amssymb}
\usepackage{amscd}
\usepackage{xspace}
\usepackage{verbatim}
\usepackage{latexsym}
\usepackage{amsfonts}
\usepackage{graphicx}
\usepackage{float}
\usepackage{color}

\newcommand{\mysection}[1]{
\section{#1}\setcounter{equation}{0}}
\title{\bf Separable infinity harmonic functions in cones}%%
\author{{\bf Marie-Fran\c{c}oise Bidaut-V\'eron} \\{\bf Marta Garcia-Huidobro}\\
 {\bf Laurent V\'eron}\\[2mm]
}%%ï¿½ï¿½Ä

\date{}
\begin{document}
 \maketitle
% \noindent{\small {\bf Abstract} We study the existence and uniqueness of  solutions of $\partial_tu-\Delta u+u^q=0$ ($q>1$) in $\Omega\times (0,\infty)$ where $\Omega\subset\mathbb R^N$ is a domain with a compact boundary, subject to the conditions $u=f\geq 0$ on $\partial\Omega\times (0,\infty)$ and the initial condition $\lim_{t\to 0}u(x,t)=\infty$. By means of Brezis' theory of maximal monotone operators in Hilbert spaces, we construct a minimal solution when $f=0$, whatever is the regularity of the boundary of the domain. When $\partial\Omega$ satisfies the parabolic Wiener criterion and $f$ is continuous, we construct a maximal solution and prove that it is the unique solution which blows-up at $t=0$.
% }

% \noindent
% {\it \footnotesize 1991 Mathematics Subject Classification}. {\scriptsize
% 35K60}.\\
% {\it \footnotesize Key words}. {\scriptsize Parabolic equations, singular solutions, semi-groups of contractions, maximal monotone operators, Wiener criterion.}
% \vspace{1mm}
% \hspace{.05in}

%% FONT commands

%% FONT commands
\newcommand{\txt}[1]{\;\text{ #1 }\;}%% Used in math only
\newcommand{\tbf}{\textbf}%% Bold face. Usage: \tbf{...}
\newcommand{\tit}{\textit}%% Italic
\newcommand{\tsc}{\textsc}%% Small caps
\newcommand{\trm}{\textrm}
\newcommand{\mbf}{\mathbf}%% Math bold
\newcommand{\mrm}{\mathrm}%% Math Roman
\newcommand{\bsym}{\boldsymbol}%% Bold math symbol
%%Macros for changing font size in math.
\newcommand{\scs}{\scriptstyle}%% as in subscript
\newcommand{\sss}{\scriptscriptstyle}%% as in sub-subscript
\newcommand{\txts}{\textstyle}
\newcommand{\dsps}{\displaystyle}
%%Macros for changing font size in text.
\newcommand{\fnz}{\footnotesize}
\newcommand{\scz}{\scriptsize}
%%\tiny<\scz<\fsz<\small<\large<\Large<\huge<\Huge
%%%%%%%%%%%%
%%%%%%%%%%%%
%% EQUATION commands
\newcommand{\be}{
\begin{equation}
}
\newcommand{\bel}[1]{
\begin{equation}
\label{#1}}
\newcommand{\ee}{
\end{equation}
}%% This macro does not work with amstex.
\newcommand{\eqnl}[2]{
\begin{equation}
\label{#1}{#2}
\end{equation}
}%%use not advisable; confusing
%%%%%%%%%%%%%%%
%% Unnumbered THEOREM env.
%% New env. to be used for unnumbered theorem, lemma etc. (but with specified name)
\newtheorem{subn}{\name}
\newcommand{\bsn}[1]{\def\name{#1}
\begin{subn}}
\newcommand{\esn}{
\end{subn}}
%%%%%%%%%%%%%%
%% NUMBERED THEOREM env.
%% Environments: theorem, lemma, corollary defintion and related commands,
%% designed to provide consecutive numbering of these forms.
\newtheorem{sub}{\name}[section]
\newcommand{\dn}[1]{\def\name{#1}}   %used in conjuction with sub or subn.
\newcommand{\bs}{
\begin{sub}}
\newcommand{\es}{
\end{sub}}
\newcommand{\bsl}[1]{
\begin{sub}\label{#1}}
%% the above must be preceeded by \dn (name definition),
%% however this is superceded by the list of commands bth etc.  below.
%%%%%%%%%%%%
%% NUMBERED THEOREM env. (cont.)
%% List of commands derived from 'sub' env. for theorem, lemma etc.
%% designed to provide consecutive numbering of these forms.
\newcommand{\bth}[1]{\def\name{Theorem}
\begin{sub}\label{t:#1}}
\newcommand{\blemma}[1]{\def\name{Lemma}
\begin{sub}\label{l:#1}}
\newcommand{\bcor}[1]{\def\name{Corollary}
\begin{sub}\label{c:#1}}
\newcommand{\bdef}[1]{\def\name{Definition}
\begin{sub}\label{d:#1}}
\newcommand{\bprop}[1]{\def\name{Proposition}
\begin{sub}\label{p:#1}}
%%%%%%%%%%%%%%%%%%%%%%%%%%%%%%%%%%
%% RERERENCE commands.
%% \newcommand{\R}[1]{$(\ref{#1})}
\newcommand{\R}{\eqref}
\newcommand{\rth}[1]{Theorem~\ref{t:#1}}
\newcommand{\rlemma}[1]{Lemma~\ref{l:#1}}
\newcommand{\rcor}[1]{Corollary~\ref{c:#1}}
\newcommand{\rdef}[1]{Definition~\ref{d:#1}}
\newcommand{\rprop}[1]{Proposition~\ref{p:#1}} 
%%%%%%%%%%%
%% ARRAY commands.
\newcommand{\BA}{
\begin{array}}
\newcommand{\EA}{
\end{array}}
\newcommand{\BAN}{\renewcommand{\arraystretch}{1.2}
\setlength{\arraycolsep}{2pt}
\begin{array}}
\newcommand{\BAV}[2]{\renewcommand{\arraystretch}{#1}
\setlength{\arraycolsep}{#2}
\begin{array}}
%Note: The first variable gives the amount of stretching: (#1) x default.
%For instance #1=1.2 means a 20% stretching. The second variable should be
%written for instance in the form  4pt ; here the default is 5pt
%\newcommand{\EAN}{\end{array}\setlength{\arraycolsep}{5pt}}
\newcommand{\BSA}{
\begin{subarray}}
\newcommand{\ESA}{\end{subarray}}
%Note: These are used in subscripts as well as superscripts. They work essentially
%% like 'array'.
\newcommand{\BAL}{\begin{aligned}}
\newcommand{\EAL}{\end{aligned}}
\newcommand{\BALG}{\begin{alignat}}
\newcommand{\EALG}{\end{alignat}}%% the abbrev. does not work with latex2e
\newcommand{\BALGN}{\begin{alignat*}}
\newcommand{\EALGN}{\end{alignat*}}%% the abbrev. does not work with latex2e
%% The 'aligned' environment must be placed inside an 'equation' env.
%% in the same way as the array.
%% One could use also the 'align' env. or the 'alignat' env.
%% However in this case each line is numbered, unless '\notag' is used.
%% The 'alignat'
%% has a slightly different format (the number of columns must be specified in advance)
%% but it has the advantage that the distance between columns is at our disposition.
%% (The default would be zero distance.) Using 'alignat*' we can have the advantages
%% of alignat plus the situation where separate lines are not numbered.
%% However in this case there is no numbering at all (unless we provide a tag).
%%%%%%%%%%
%% PROOF, REMARK etc.
\newcommand{\note}[1]{\textit{#1.}\hspace{2mm}}
\newcommand{\Proof}{\note{Proof}}
\newcommand{\qeda}{\hspace{10mm}\hfill $\square$}
\newcommand{\Remark}{\note{Remark}}
%%%%%%%% Style command.
\newcommand{\modin}{$\,$\\
[-4mm] \indent}
%% To be used after \mysection in order to start new line with \indent.
%%%%%%%%%%%%
%% MATHEMATICAL symbols
\newcommand{\forevery}{\quad \forall}
\newcommand{\set}[1]{\{#1\}}
\newcommand{\setdef}[2]{\{\,#1:\,#2\,\}}
\newcommand{\setm}[2]{\{\,#1\mid #2\,\}}
%% Arrows
\newcommand{\lra}{\longrightarrow}
\newcommand{\sgn}{\rm{sgn}}
\newcommand{\lla}{\longleftarrow}
\newcommand{\llra}{\longleftrightarrow}
\newcommand{\Lra}{\Longrightarrow}
\newcommand{\Lla}{\Longleftarrow}
\newcommand{\Llra}{\Longleftrightarrow}
\newcommand{\warrow}{\rightharpoonup}
%% Brackets, delimiters
\newcommand{
\paran}[1]{\left (#1 \right )}%% adjustable parantheses
\newcommand{\sqbr}[1]{\left [#1 \right ]}%% adjustable square brackets
\newcommand{\curlybr}[1]{\left \{#1 \right \}}%% adjustable curly brackets
\newcommand{\abs}[1]{\left |#1\right |}%% adjustable vertical delimiters
\newcommand{\norm}[1]{\left \|#1\right \|}%% adjustable norm
\newcommand{
\paranb}[1]{\big (#1 \big )}%% non-adjustable parantheses (big)
\newcommand{\lsqbrb}[1]{\big [#1 \big ]}%% non-adjustable square brackets (big)
\newcommand{\lcurlybrb}[1]{\big \{#1 \big \}}%% non-adjustable curly brackets (big)
\newcommand{\absb}[1]{\big |#1\big |}%% non-adjustable vertical delimiters (big)
\newcommand{\normb}[1]{\big \|#1\big \|}%% non-adjustable norm (big)
\newcommand{
\paranB}[1]{\Big (#1 \Big )}%% non-adjustable parantheses (Big)
\newcommand{\absB}[1]{\Big |#1\Big |}%% non-adjustable vertical delimiters (Big)
\newcommand{\normB}[1]{\Big \|#1\Big \|}%% non-adjustable norm (Big)

%%%%%%%%%%%%%%%%%
%% Adjustable parantheses etc. in a different DEFINITION format.
%\def\adp(#1){\left (#1 \right )}%% adjustable parantheses
%\def\adsb(#1){\left [#1\right ]}%% adjustable square brackets
%\def\adcb(#1){\left \{#1\right \}}%% adjustable curly brackets
%\def\abs|#1|{\left |#1\right |}%% adjustable vertical delimiters
%%%%%%%%%%%%%%%%
%% More mathematical symbols
\newcommand{\thkl}{\rule[-.5mm]{.3mm}{3mm}}
\newcommand{\thknorm}[1]{\thkl #1 \thkl\,}
\newcommand{\trinorm}[1]{|\!|\!| #1 |\!|\!|\,}
\newcommand{\bang}[1]{\langle #1 \rangle}%% angle bracket
\def\angb<#1>{\langle #1 \rangle}%% angle bracket
%% The two last lines yield the same result.
%% The second is used as follows: \angb<a,b>
\newcommand{\vstrut}[1]{\rule{0mm}{#1}}
\newcommand{\rec}[1]{\frac{1}{#1}}
%% OPERATOR names.
%% OPERATOR names.
\newcommand{\opname}[1]{\mbox{\rm #1}\,}
\newcommand{\supp}{\opname{supp}}
\newcommand{\dist}{\opname{dist}}
\newcommand{\myfrac}[2]{{\displaystyle \frac{#1}{#2} }}
\newcommand{\myint}[2]{{\displaystyle \int_{#1}^{#2}}}
\newcommand{\mysum}[2]{{\displaystyle \sum_{#1}^{#2}}}
\newcommand {\dint}{{\displaystyle \int\!\!\int}}%%%%%%%%%%
%%%%%%% SPACE commands
\newcommand{\q}{\quad}
\newcommand{\qq}{\qquad}
\newcommand{\hsp}[1]{\hspace{#1mm}}
\newcommand{\vsp}[1]{\vspace{#1mm}}
%%%%%%%%%%%
%% ABREVIATIONS
\newcommand{\ity}{\infty}
\newcommand{\prt}{
\partial}
\newcommand{\sms}{\setminus}
\newcommand{\ems}{\emptyset}
\newcommand{\ti}{\times}
\newcommand{\pr}{^\prime}
\newcommand{\ppr}{^{\prime\prime}}
\newcommand{\tl}{\tilde}
\newcommand{\sbs}{\subset}
\newcommand{\sbeq}{\subseteq}
\newcommand{\nind}{\noindent}
\newcommand{\ind}{\indent}
\newcommand{\ovl}{\overline}
\newcommand{\unl}{\underline}
\newcommand{\nin}{\not\in}
\newcommand{\pfrac}[2]{\genfrac{(}{)}{}{}{#1}{#2}}% frac with parantheses.
%%%%%%%%%%%
%%%%%%%%%%%%%

%%Macros for Greek letters.
\def\ga{\alpha}     \def\gb{\beta}       \def\gg{\gamma}
\def\gc{\chi}       \def\gd{\delta}      \def\ge{\epsilon}
\def\gth{\theta}                         \def\vge{\varepsilon}
\def\gf{\phi}       \def\vgf{\varphi}    \def\gh{\eta}
\def\gi{\iota}      \def\gk{\kappa}      \def\gl{\lambda}
\def\gm{\mu}        \def\gn{\nu}         \def\gp{\pi}
\def\vgp{\varpi}    \def\gr{\rho}        \def\vgr{\varrho}
\def\gs{\sigma}     \def\vgs{\varsigma}  \def\gt{\tau}
\def\gu{\upsilon}   \def\gv{\vartheta}   \def\gw{\omega}
\def\gx{\xi}        \def\gy{\psi}        \def\gz{\zeta}
\def\Gg{\Gamma}     \def\Gd{\Delta}      \def\Gf{\Phi}
\def\Gth{\Theta}
\def\Gl{\Lambda}    \def\Gs{\Sigma}      \def\Gp{\Pi}
\def\Gw{\Omega}     \def\Gx{\Xi}         \def\Gy{\Psi}

%%Macros for calligraphic letters.
\def\CS{{\mathcal S}}   \def\CM{{\mathcal M}}   \def\CN{{\mathcal N}}
\def\CR{{\mathcal R}}   \def\CO{{\mathcal O}}   \def\CP{{\mathcal P}}
\def\CA{{\mathcal A}}   \def\CB{{\mathcal B}}   \def\CC{{\mathcal C}}
\def\CD{{\mathcal D}}   \def\CE{{\mathcal E}}   \def\CF{{\mathcal F}}
\def\CG{{\mathcal G}}   \def\CH{{\mathcal H}}   \def\CI{{\mathcal I}}
\def\CJ{{\mathcal J}}   \def\CK{{\mathcal K}}   \def\CL{{\mathcal L}}
\def\CT{{\mathcal T}}   \def\CU{{\mathcal U}}   \def\CV{{\mathcal V}}
\def\CZ{{\mathcal Z}}   \def\CX{{\mathcal X}}   \def\CY{{\mathcal Y}}
\def\CW{{\mathcal W}} \def\CQ{{\mathcal Q}} 
%%%%%
%%Macros for 'blackboard' letters (See (27) for display.)
\def\BBA {\mathbb A}   \def\BBb {\mathbb B}    \def\BBC {\mathbb C}
\def\BBD {\mathbb D}   \def\BBE {\mathbb E}    \def\BBF {\mathbb F}
\def\BBG {\mathbb G}   \def\BBH {\mathbb H}    \def\BBI {\mathbb I}
\def\BBJ {\mathbb J}   \def\BBK {\mathbb K}    \def\BBL {\mathbb L}
\def\BBM {\mathbb M}   \def\BBN {\mathbb N}    \def\BBO {\mathbb O}
\def\BBP {\mathbb P}   \def\BBR {\mathbb R}    \def\BBS {\mathbb S}
\def\BBT {\mathbb T}   \def\BBU {\mathbb U}    \def\BBV {\mathbb V}
\def\BBW {\mathbb W}   \def\BBX {\mathbb X}    \def\BBY {\mathbb Y}
\def\BBZ {\mathbb Z}

%%Macros for Ghotic (Fraktur) letters.
\def\GTA {\mathfrak A}   \def\GTB {\mathfrak B}    \def\GTC {\mathfrak C}
\def\GTD {\mathfrak D}   \def\GTE {\mathfrak E}    \def\GTF {\mathfrak F}
\def\GTG {\mathfrak G}   \def\GTH {\mathfrak H}    \def\GTI {\mathfrak I}
\def\GTJ {\mathfrak J}   \def\GTK {\mathfrak K}    \def\GTL {\mathfrak L}
\def\GTM {\mathfrak M}   \def\GTN {\mathfrak N}    \def\GTO {\mathfrak O}
\def\GTP {\mathfrak P}   \def\GTR {\mathfrak R}    \def\GTS {\mathfrak S}
\def\GTT {\mathfrak T}   \def\GTU {\mathfrak U}    \def\GTV {\mathfrak V}
\def\GTW {\mathfrak W}   \def\GTX {\mathfrak X}    \def\GTY {\mathfrak Y}
\def\GTZ {\mathfrak Z}   \def\GTQ {\mathfrak Q}

\font\Sym= msam10 % special symbols
\def\SYM#1{\hbox{\Sym #1}}
\newcommand{\bdw}{\prt\Gw\xspace}
\medskip
\noindent{\small {\bf Abstract}
We study the existence of separable infinity harmonic functions in any cone of $\BBR^N$ vanishing on its boundary under the form $u(r,\gs)=r^{-\gb}\psi(\gs)$. We prove that such solutions exist, the spherical part $\psi$ satisfies a nonlinear eigenvalue problem on a subdomain of the sphere $S^{N-1}$ and that the  exponents $\gb=\gb_+>0$ and $\gb=\gb_-<0$ are uniquely determined if the domain is smooth.  We extend some of our results to non-smooth domains.
}\smallskip

\tableofcontents
\vspace{1mm}
%%%%%%%%%%%%%%%%%%%%%%%%%%%%%%%%%%%%%%%%%%%%%%%%%%%%%%

%%%%%%%%%%%%%%%%%%%%%%%%%%%%%%%%%%%%%%%%%%%%%%%%%%%%%%%%%%%%%%%%%%%%%%%%%%%%%%%%%%%%%%%%%%%%%%%%%%%%%%%%%%%%%%%%%%%%%%%SECTION-THE REGULAR CASE%%%%%%%%%%%%%%%%%%%%%%%%%%%%%%%%%%%%%%%%%%%%%%%%%%%%%%%%%%%%%%%%%%%%%%%%%%%%%%%%%%%%%%%%%%%%%%%%%%%%%%%%%%%%%%%%%%%%%%%%%%%%%%%%%%%%%%

\mysection{Introduction}

Let $S$ be a $C^3$ subdomain of the unit sphere $S^{N-1}$ of $\BBR^N$ and $C_S:=\{\gl\gs\in\BBR^N:\gl>0,\,\gs\in S\}$ is the positive cone generated by $S$. In this paper we study the existence of positive solutions of 
\bel{I1}
\Gd_\infty u:=\myfrac{1}{2}\nabla \abs{\nabla u}^2.\nabla u=0
\ee
in $C_S$  vanishing on $\prt C_S\setminus\{0\}$ under the form 
\bel{I2}
u(x)=u(r,\gs)=r^{-\gb}\psi(\gs),
\ee
where $\gb\in\BBR$ and $(r,\gs)\in \BBR_+\ti S^{N-1}$ are the spherical coordinates $\BBR^N$; such a function $u$ is called a {\it separable infinity harmonic function}. The function $\psi$ satisfies the {\it spherical infinity harmonic problem in $S$}  
\bel{I3}\BA {ll}
\myfrac{1}{2}\nabla'\abs{\nabla'\psi}^2.\nabla'\psi+\gb(2\gb+1)\abs{\nabla'\psi}^2\psi+\gb^3(\gb+1)\psi^3=0&\text{in } S\\
\phantom{\myfrac{1}{2}\nabla'\abs{\nabla'\psi}^2.\nabla'\psi+\gb(2\gb+1)\abs{\nabla'\psi}^2\psi+\gb^3(\gb+1)^3}
\psi=0\qquad&\text{on }\prt S,
\EA\ee
where $\nabla'$ is the covariant gradient on $S^{N-1}$ for the canonical metric and $(a,b)\mapsto a.b$ the associated quadratic form. The role of the infinity Laplacian for Lipschitz extension of Lipchitz continuous functions defined in a domain has been pointed out by Aronsson in his seminal paper \cite{Aro}. When the infinity Laplacian $\Gd_\infty$ is replaced by the $p$-Laplacian, the research of regular ($\gb<0$)  separable $p$-harmonic functions has been carried out by Krol \cite{Kro} in the $2$-dim case and by Tolksdorff \cite{Tolk} in the general case. Following Krol' s method, Kichenassamy and V\'eron \cite {KiVe} studied the $2$-dim singular case ($\gb>0$). Finally, by a completely different approach and in a more general setting Porretta and V\'eron \cite{PoVe} studied the general case. In that case, the function $\psi$ 
satisfies the {\it spherical $p$-harmonic problem in $S$}
\bel{I4}\BA {ll}
div'\left((\gb^2\psi^2+|\nabla' \psi|^2)^{\frac{p-2}{2}}\nabla' \psi\right)+\gb\gl_\gb(\gb^2\psi^2+|\nabla' \psi|^2)^{\frac{p-2}{2}} \psi=0&\text{in } S\\
\phantom{div'\left((\gb^2\psi^2+|\nabla' \psi|^2)^{\frac{p-2}{2}}\nabla' \psi\right)+\gb\gl_\gb(\gb^2\psi^2+|\nabla' \psi|^2)^{\frac{p-2}{2}} }
\psi=0\quad&\text{on }\prt S,
\EA\ee
where $\gl_\gb=\gb(p-1)+p-N$ and $div'$ is the divergence operator acting on vector fields in $TS^{N-1}$. \smallskip

Following an idea which was introduced by Lasry and Lions \cite{LaLi}, Porretta-V\'eron's method was to transform the equation $(\ref{I4})$ by setting
\bel{I5}\BA {ll}
w=-\myfrac{1}{\gb}\ln \psi
\EA\ee
in the case $\gb>0$. The function $w$ satisfies the new problem
\bel{I6}\BA {ll}\displaystyle
 -div'\left(\left(1+|\nabla' w|^2\right)^{p/2-1}\nabla' w\right)+\left(1+|\nabla' w|^2\right)^{p/2-1}\left(\gb(p-1)|\nabla' w|^2+\gl_\gb\right)=0\quad\mbox{ in }S
 \\[2mm]\phantom{-----------------,---------,-}\displaystyle
 \lim_{ \gr(\gs)\to 0}w( \gs)=\infty,
\EA\ee
where $\gr(\gs):=\dist(\gs,\prt S)$ is the distance understood in the the geodesic sense on $S$.\\

In this article we borrow ideas used in \cite{PoVe} to transform problem $(\ref{I1})$ by introducing the function $w$ defined by $(\ref{I5})$. Then $w$ satisfies, in the viscosity sense,
\bel{I7}\BA {ll}\displaystyle
-\myfrac{1}{2}\nabla'\abs{\nabla'w}^2.\nabla'w+\gb\abs{\nabla'w}^4+(2\gb+1)\abs{\nabla'w}^2
+\gb+1=0\quad\mbox{ in }S
 \\[2mm]\phantom{-----------,----+\gb+1-}\displaystyle
 \lim_{ \gr(\gs)\to 0}w( \gs)=\infty.
\EA\ee
We first prove\smallskip

\nind{\bf Theorem A}. {\it Let $S\subset S^{N-1}$ be a proper subdomain of $S^{N-1}$ with a $C^3$ boundary. Then for any $\gb>0$ there exists a locally Lipschitz continuous function $w$ and a unique $\gl(\gb)>0$ satisfying in the viscosity sense
\bel{I7+}\BA {ll}\displaystyle
-\myfrac{1}{2}\nabla'\abs{\nabla'w}^2.\nabla'w+\gb\abs{\nabla'w}^4+(2\gb+1)\abs{\nabla'w}^2+\gl(\gb)
=0\quad\mbox{ in }S
 \\[2mm]\phantom{---------+\gl(\gb)-------,}\displaystyle
 \lim_{ \gr(\gs)\to 0}w( \gs)=\infty.
\EA\ee
where $\gr(.)$ is the geodesic distance from points in $S$ to $\prt S$. 
}\smallskip

Then we prove that there exists a unique $\gb$ such that $\gl(\gb)=\gb+1$. In a similar way we study the regular case where $\gb<0$ in $(\ref{I2})$, (we denote $-\gb=\gm>0$), and we obtain\smallskip

\nind{\bf Theorem B}. {\it Let $S\subset S^{N-1}$ be subdomain with a $C^3$ boundary. Then there exist exactly two real positive numbers $\gb_{s}$ and $\gm_s$ and at least two positive functions $\psi_s$ and 
$\gw_s$ in $C(\overline S)\cap C_{loc}^{0,1}(S)$ (up to multiplication by constants) such that the two functions $u_{s,{\scriptscriptstyle+}}$ and 
$u_{s,{\scriptscriptstyle-}}$ defined in $C_S$ by $u_{s,{\scriptscriptstyle+}}(r,\gs):=r^{-\gb_s}\psi_s(\gs)$ and $u_{s,{\scriptscriptstyle-}}(r,\gs):=r^{\gm_s}\gw_s(\gs)$ are infinity harmonic in $C_S$ and vanish on $\prt C_S\setminus\{0\}$ and $\prt C_S$ respectively. Furthermore $\gb_s$ and $\gm_s$ are decreasing functions of $S$ for the inclusion order relation on sets.
}\smallskip

The previous results can be extended to general regular domains on a Riemannian manifold as in \cite{PoVe}. It is an open problem whether the positive solutions associated to the same exponent $\gb_s$ (or $b_s$) are proportional, see discussion in Remark p. 15.

\smallskip

In the special case of a rotationally symmetric domain $S$ we have a more precise result which allows us to characterize all the separable infinity harmonic functions in $C_S$ which keep a constant sign and vanish on $\prt C_S\setminus\{0\}$. We denote by $\phi\in (0,\gp)$ the azimuthal angle from the North pole $N$ on $S^{N-1}$.\smallskip

\nind{\bf Theorem C}. {\it Let $S_\ga$ be the spherical cap with azimuthal opening $\ga\in (0,\gp]$. Then there exist two positive functions $\psi_\ga$ and $\gw_\ga$ in $C^{\infty}(\overline S)$, vanishing on $\prt S$, such that the two functions
\bel{I8}\BA {ll}\displaystyle
u_{\ga,+}(r,\gs)=r^{-\frac{\gp^2}{4\ga(\gp+\ga)}}\psi_\ga(\gs),
\EA\ee
and
\bel{I8-1}\BA {ll}\displaystyle
u_{\ga,-}(r,\gs)=r^{\frac{\gp^2}{4\ga(\gp-\ga)}}\gw_\ga(\gs),
\EA\ee
are infinity harmonic in $C_{S_\ga}$ and vanish on $\prt C_{S_\ga}\setminus\{0\}$. The two functions $\psi_\ga$ and $\gw_\ga$ depend only on the variable $\gf\in (0,\ga]$ and are unique in the class of rotanionnaly invariant solutions up to multiplication by constants.}
\smallskip

This study reduced to an ordinary differential equation which has been already treated by T. Bhattacharya in \cite{Bat} and \cite{Bat1}. But for the sake of completeness and for some related problems we present it in Section 3 of the present paper. \smallskip

Using these previous results we prove the existence of separable infinity harmonic functions in any cone $C_S$.\smallskip

\nind{\bf Theorem D}. {\it Assume $S\subsetneq S^{N-1}$ is any domain. Then there exist $\overline\gb_s>0$ and a positive function $\overline \psi_s$ in $C(\overline S)$, locally Lipschitz continuous in $S$ and vanishing on $\prt S$, such that the function 
\bel{I9}\BA {ll}\displaystyle
\overline u_{s,{\scriptscriptstyle+}}(r,\gs)=r^{-\overline\gb_s}\overline \psi_s(\gs),
\EA\ee
is infinity harmonic in $C_{S}$ and vanishes on $\prt C_{S}\setminus\{0\}$.} \smallskip

When the cone $C_S$  is a little more regular, the construction of the spherical infinity harmonic functions can be performed via an approximation from outside.\smallskip

\nind{\bf Theorem E}. {\it Assume that $S\subsetneq S^{N-1}$ is an outward accessible domain, i.e. $\prt S=\prt\overline S^c$. Then there exist $\underline\gb_s\in(0,\overline\gb_s]$ and a positive function $\underline \psi_s\in C(\overline S)$, locally Lipschitz continuous in $S$ and vanishing on $\prt S$, such that the function 
\bel{I9-1}\BA {ll}\displaystyle
\underline u_{s,{\scriptscriptstyle+}}(r,\gs)=r^{-\underline\gb_s}\underline \psi_s(\gs)
\EA\ee
is infinity harmonic in $C_{S}$, vanishes on $\prt C_{S}\setminus\{0\}$ and has the property that for any separable infinity-harmonic function in $C_S$ under the form $u(r,\gs)=r^{-\gb}\psi(\gs)$ where $\gb>0$ and $\psi$ in $C(\overline S)$ vanishes on $\prt S$, there holds ${\underline\gb}_s\leq \gb\leq \overline\gb_s$. }\smallskip

The uniqueness of the exponent $\gb_s$ is proved under Lipschitz  and geometric conditions on $S$.\smallskip

\nind{\bf Theorem F}. {\it Assume that $S\subsetneq S^{N-1}$ is a Lipschitz domain satifying the interior sphere condition. Then $\overline\gb_s=\underline\gb_s$. Furthermore there exists a constant $c=c(S,p)>0$ such that for any two positive functions $\psi_i$, $i=1,2$, satisfying the spherical p-harmonic problem $(\ref{I4})$, there holds
\bel{I10}\BA {ll}\displaystyle
\myfrac{\psi_1(\gs)}{\psi_2(\gs)}\leq c\myfrac{\psi_1(\gs')}{\psi_2(\gs')}\qquad\forall (\gs,\gs')\in S.
\EA\ee
}\medskip

Note that  the statements and the proofs of Theorems D, E and F can be easily modified if one considers regular infinity harmonic functions in $C_{S}$ which vanish on $\prt C_{S}$.
\smallskip

\nind{\bf Acknowledgements.} This article has been prepared with the support of the collaboration
programs ECOS C14E08 and FONDECYT grant 1160540 for the three authors. The authors are grateful to the referee for a careful reading of their work.
\bigskip

\mysection{The smooth case}
\medskip

We assume that $(r,\gs)\in\BBR_+\ti S^{N-1}$ are the spherical coordinates of $x\in\BBR^{N}$. If $u$ is a $C^1$ function, then $\nabla u=u_r{\bf e}+\frac{1}{r}\nabla'u$ where ${\bf e}=\frac x{|x|}$ and $\nabla'$ is the tangential gradient of $u(r,.)$ identified to the covariant gradient thanks to the canonical imbedding of $S^{N-1}$ into $\BBR^N$. Then $|\nabla u|^2=u_r^2+\frac{1}{r^2}|\nabla'u|^2$, thus
$$-\Gd_\infty u=
-\left(u_r^2+\frac{1}{r^2}|\nabla'u|^2\right)_ru_r-\frac{1}{r^2}\nabla '\left(u_r^2+\frac{1}{r^2}|\nabla'u|^2\right).\nabla'u=0.
$$
A solution $-\Gd_\infty u=0$ which has the form $u(x)=u(r,\gs)=r^{-\gb}\psi(\gs)$ satisfies, in the viscosity sense, the spherical infinity harmonic equation
\bel{G1}\BA {ll}
\myfrac{1}{2}\nabla'\abs{\nabla'\psi}^2.\nabla'\psi+\gb(2\gb+1)\abs{\nabla'\psi}^2\psi+\gb^3(\gb+1)\psi^3=0.
\EA
\ee

\bth{main} For any $C^3$ domain $S\subset S^{N-1}$ there exists a unique $\gb_s>0$
and one nonnegative function $\psi\in C^{0,1}(\overline S)$ solution of 
\bel{E5}\BA {ll}
-\myfrac{1}{2}\nabla'|\nabla' \psi|^2.\nabla' \psi=\gb(2\gb+1)|\nabla' \psi|^2w+\gb^3(\gb+1)\psi^3\quad&\text{in }\; S
\\\phantom{-\myfrac{1}{2}\nabla'|\nabla' \psi|^2.\nabla' }
\psi=0&\text{in }\; \prt S. 
\EA\ee
such that  the function $(r,\gs)\mapsto u_s(r,\gs):= r^{-\gb_s}\psi(\gs)$ is positive and $\infty$-harmonic    in the cone $C_s=\{x=\gl\gs\in \BBR^N:\gl>0,\, \gs\in S\}$ and vanish on  $\prt S\setminus\{0\}$.\es

Following Porretta-Veron's method, we transform the eigenvalue problem into a large solution problem with absorption by setting
\bel{E5-1}\BA {ll}w=-\myfrac{1}{\gb}\ln \psi.
\EA\ee
Therefore the formal new problem is to prove the existence of a unique $\gb>0$ and of a nonnegative function $w$ such that
\bel{E6}\BA {ll}
\myfrac{1}{2}\nabla'|\nabla' w|^2.\nabla' w-\gb|\nabla' w|^4-(2\gb+1)|\nabla' w|^2=\gb+1\quad&\text{in }\; S
\\\phantom{\myfrac{1}{2}\nabla'|\nabla'|^2.\nabla' w-\gb|\nabla' w|^4-(2\gb+1)|\nabla' w|^2}
w=\infty&\text{in }\; \prt S. 
\EA\ee
The two problems are clearly equivalent for $C^2$ solutions. Since the mapping $w\mapsto \psi$ is smooth and decreasing, it exchanges supersolutions (resp. subsolutions) into subsolutions (resp. supersolutions). Therefore the two problems $(\ref{E5})$-$(\ref{E6})$ are also equivalent if we deal with continuous viscosity solutions. \smallskip

In order to increase the regularity of the solutions and to avoid the difficulties coming form the fact the above problem is invariant if we add a constant to a solution, instead of $(\ref{E6})$ we consider the regularized problem with absorption

\bel{E7}\BA {ll}
-\gd\Gd w-\myfrac{1}{2}\nabla'|\nabla' w|^2.\nabla' w+\gg|\nabla' w|^4+(2\gg+1)|\nabla' w|^2+\ge w=0\quad&\text{in }\; S
\\\phantom{-\gd\Gd w-\myfrac{1}{2}\nabla'|\nabla' w|^2.\nabla' w-\gg|\nabla' w|^4-(2\gg+1)|\nabla' w|^2-\ge}
w=\infty&\text{in }\; \prt S,
\EA\ee
where $\ge,\gd$ are two positive parameters. We will obtain below local estimates on $\nabla' w$ independent of $\ge$ and $\gd$. Thanks to these estimates we will let successively $\gd$ and $\ge$ to $0$ and obtain that, up to a constant, the term $\ge w$ converges to some unique $\gl(\gg)$ called {\it the  ergodic constant} although it has a probabilistic interpretation only in the case of the ordinary Laplacian \cite{LaLi}. The limit problem of $(\ref{E7})$ is  the following

\bel{E8}\BA {ll}
-\myfrac{1}{2}\nabla'|\nabla' w|^2.\nabla' w+\gg|\nabla' w|^4+(2\gg+1)|\nabla' w|^2+\gl(\gg)=0\quad&\text{in }\; S
\\\phantom{-\myfrac{1}{2}\nabla'|\nabla' w|^2.\nabla' +\gg|\nabla' w|^4+(2\gg+1)|\nabla' w|^2+\gl(\gg)}
w=\infty&\text{in }\; \prt S.
\EA\ee

\subsection{Two-sided estimates}

We denote the "positive" geodesic distance $\gr(\gs)=\dist (\gs,\prt S)$. If $\gs\in S^c$ we set $\tilde\gr(\gs)=-\dist (\gs,\prt S)$. If $\gs_1$ and $\gs_2$ are not antipodal points there exists a unique minimizing geodesic between $\gs_1$ and $\gs_2$. It is an arc of a Riemannian circle (or great circle). The geodesic distance between 
$\gs_1$ and $\gs_2$ is denoted by $\ell (\gs_1,\gs_2)$. It coincides with the angle determined by the two straight lines from $0$ to $\gs_1$ and $0$ to $\gs_2$.  
At this point it is convenient to use Fermi coordinates in $S$ in a neighborhood of $\prt S$. We set 
$$S_\gt=\{\gs\in S:\gr(\gs)<\gt\}\,,\;S'_\gt=S\setminus\overline S_\gt\,,\;\Gs_\gt= \{\gs\in S:\gr(\gs)=\gt\}.$$
If $\gt\leq\gt_0$ for any $\gs\in S_\gt$ there exists a unique $z_\gs\in \prt S$ such that $\ell (\gs,z_\gs)=\gr(\gs)$.
These Fermi coordinates of $\gs$ are defined by $(\gt,z)\in [0,\gt_0)\ti \prt S$. The mapping  $\Gp$ such that
$$\Gp(\gs)=(\gr(\gs),z_\gs)\qquad\forall \gs\in S_{\gt_0}, 
$$
is a $C^2$ diffeomorphism from $S_{\gt_0}$ into $[0,\gt_0)\ti \prt S$. The expression of the Laplace-Beltrami operator in $S_{\gt_0}$ is given in \cite{BanMar}:
\bel{E10-1}\BA {ll}
\Gd'u(\gs)=\myfrac{\prt^2u}{\prt\gt^2}-(N-2)H\myfrac{\prt u}{\prt\gt}+\tilde\Gd'_z u\qquad\forall \gs=\Gp^{-1}((\gt,z)),
\EA\ee
where $H=H(\gt,z)$ is the mean curvature of $\Gs_\gt$ and $\tilde\Gd'_z$ is a second order elliptic operator acting on functions defined on $\Gs_\gt$. If $g=(g_{ij})$ is the metric tensor on $S^{N-1}$ and by convention, $|g|=det(g_{ij})$, this operator  admits the following expression 
$$\tilde\Gd'_zu=\myfrac{1}{\sqrt {|g|}}\sum_{j=1}^{N-2}\myfrac{\prt}{\prt z_j}\left(\sqrt {|g| }a_{j}\myfrac{\prt u}{\prt z_j}\right),
$$
for some $a_j>0$ if we take for coordinates curve-frame $z_j$ a system of orthogonal 1-dim great circles on $\Gg$ intersecting at $z_\gs$ (these circle corresponds to the ($N-2$)-principal curvatures at this points). The coefficients $a_j$ depend both on $z$ and $\gt$. Thus, if $u$ depends only on $\gr$, 
\bel{E10-2}\BA {ll}
\Gd'u(\gs)=\myfrac{\prt^2u}{\prt\gt^2}-(N-2)H\myfrac{\prt u}{\prt\gt}.
\EA\ee
The expression of $H$ is given in \cite{BanMar} and we can assume that $\gt_0$ is small enough so that $H$ remains bounded.
\medskip

We extend the geodesic distance $\gr(x)=\dist (x,\prt S)$ as a smooth positive function so that $\tilde \rho(x):=\rho(x)$ if $\gr(x)\leq \gt_0$ and thus, it the same neighborhood of $\prt S$, $\nabla\tilde \rho(x)={\bf n}_{z_x}$, the unit outward normal vector to $\prt S$ at the point $z_x=Proj_{\prt\Gw}(x)$. \smallskip

If $w$ depends only on $\gr$, $(\ref{E7})$ becomes 
\bel{E11}\BA {ll}
\gd w''-(N-2)Hw'+w'^2w''-\gg w'^4-(2\gg+1)w'^2-\ge w=0\quad&\text{in }\; S_{\gt_0}
\\\phantom{\gt w''-(N-2)Hw'+w'^2w'-\gg|w'^4-(2\gg+1)w'^2-\ge}
w=\infty&\text{in }\; \prt S.
\EA\ee 
In the sequel we put
\bel{E10}\BA {ll}
\CP_\gd(u):=-\gd\Gd u-\myfrac{1}{2}\nabla'|\nabla' u|^2.\nabla' u+\gg|\nabla' u|^4+(2\gg+1)|\nabla' u|^2+\ge u=\tilde \CP_\gd
(u)+\ge u.
\EA\ee

%\subsection{Estimate on the function}
%%%%%%%%%%%%%%%%%%%%%%%%%%%%%%%%%%%%%%%%%%%%%%%%%%%%%%%%%%%%%%%%%%%%%%%%%%%%%%%%%%%%%%%%%%%%%%%%%%%%%%%%%%%%%%%%%%%%%%%%%%%%%%%%%%%%%%%%%%%%%%%%%%%%%%%%%%%%%%%%%%
\bprop{supsub} There exist $\gt_1\in (0,\gt_0]$, three positive constants $M$,  $\ge_0$ and $\gd_0$ and two positive functions $w^*,w_*\in C^2(S)$ such that 
$w^*>w_*$ in $S_{\gt}$, $w^*+\myfrac{1}{\gg}\ln\gr\in L^\infty(S)$ and $w_*+\myfrac{1}{\gg}\ln\gr\in L^\infty(S)$ with the property that for any $\ge\in (0,\ge_0]$ and $\gd\in (0,\gd_0]$ the two functions 

\bel{E12}\BA {ll}
\bar w(\gs)=w^*+\myfrac{M}{\ge}
\EA\ee
and 
\bel{E13}\BA {ll}
\underline w(\gs)=w_*-\myfrac{M}{\ge}
\EA\ee
are respectively a supersolution and a subsolution of $\CP_\gd(u)=0$. Furthermore any solution $w$ of problem         $(\ref{E8})$ satisfies $\underline w\leq w\leq \bar w$.
\es
\nind\Proof 
Let $a>0$. We first notice by a standard computation that the solutions of the ODE,
\bel{E14}\BA {ll}
\gd w''+w'^2w''-\gg w'^4-aw'^2=0\quad&\text{in }\; (0,1)
\\\phantom{\gd w''+w^2w''-^4-aw'^2}
w'(0)=-\infty
\EA\ee
are negative and given implicitely by 
\bel{E15}\BA {ll}
\myfrac{\gd}{a w'(\gr)}+\left(\myfrac{1}{\gg}-\myfrac{\gd}{a}\right)\sqrt{\myfrac{\gg}{a}}\tan^{-1}\left(\sqrt{\myfrac{a}{\gg}}\myfrac{1}{w'(\gr)}\right)=-\gr.
\EA\ee
In order to have a global estimate, we set $w'=-z^{-1}$, thus $(\ref{E15})$ becomes
\bel{E15-1}\BA {ll}
\myfrac{\gd z}{a }+\left(\myfrac{1}{\gg}-\myfrac{\gd}{a}\right)\sqrt{\myfrac{\gg}{a}}\tan^{-1}\left(z\sqrt{\myfrac{a}{\gg}}\right)=\gr,
\EA\ee
provided $a>\gg\gd$. Since $\tan^{-1}\left(z\sqrt{\frac{a}{\gg}}\right)\leq z\sqrt{\frac{a}{\gg}}$, we derive
\bel{E15-2}\BA {ll}
z(\gr)\geq\gg\gr\Longleftrightarrow 0>w'(\gr)\geq -\myfrac{1}{\gg \gr}\qquad\forall \gr>0,
\EA\ee
with equality only if $\gr=0$. Since we can write $(\ref{E15-1})$ as
\bel{E15-3}\BA {ll}
\sqrt{\myfrac{\gg}{a}}\myfrac{\gd }{a }\left(z\sqrt{\myfrac{a}{\gg}}\right)+\left(\myfrac{1}{\gg}-\myfrac{\gd}{a}\right)\sqrt{\myfrac{\gg}{a}}\tan^{-1}\left(z\sqrt{\myfrac{a}{\gg}}\right)=\gr\geq \tan^{-1}\left(z\sqrt{\myfrac{a}{\gg}}\right),
\EA\ee
we obtain
\bel{E15-4}\BA {ll}
z\leq \sqrt{\myfrac{\gg}{a}}\tan\left(\gr\sqrt{a\gg}\right)\Longleftrightarrow w'(\gr)\leq -\sqrt{\myfrac{a}{\gg}}\cot\left(\gr\sqrt{a\gg}\right)\qquad\forall \gr>0.
\EA\ee
Finally
\bel{E15-5}\BA {ll}
-\myfrac{1}{\gg \gr}\leq  w'(\gr)\leq -\sqrt{\myfrac{a}{\gg}}\cot\left(\gr\sqrt{a\gg}\right)\qquad\forall \gr\in (0,\gt_0].
\EA\ee
and in particular, for any $\gt_1\in (0,\gt_0]$,

\bel{E15-6}\BA {ll}
\abs{w'(\gr)}\geq \sqrt{\myfrac{a}{\gg}}\cot\left(\gt_1\sqrt{a\gg}\right)\qquad\forall \gr\in (0,\gt_1].
\EA\ee

From this estimate we derive
\bel{E15-7}\BA {ll}
w(\gr_0)+\myfrac{1}{\gg}\ln\left(\myfrac{\sin\gr_0\sqrt{a\gg}}{\sin\gr\sqrt{a\gg}}\right)\leq  w(\gr)\leq w(\gr_0)+\myfrac{1}{\gg }\ln\left(\myfrac{\gr_0}{\gr}\right)\qquad\forall \gr\in (0,\gt_1].
\EA\ee

The solution $w$ depends on the value of $a$ and $\gd$. Since $\prt S$ is smooth, we can assume that $(N-2)\abs H$ is bounded by some constant $m\geq 0$ in $S_{\gt_0}$. Denote by $w_{\gt}$ a solution satisfying $w(\gt)=0$, then it is positive in $S_{\gt}$ and 
$$\BA{lll}
\CP_\gd (w_{\gt})=\gd(N-2) H w_{\gt}'+(2\gg+1-a)w_{\gt}'^2+\ge w_{\gt}\\[2mm]
\phantom{\CP_\gd (w_{\gt})}
\geq |w_{\gt}'|((2\gg+1-a)|w_{\gt}'|-m)\\[2mm]
\phantom{\CP_\gd (w_{\gt})}
\geq |w_{\gt}'|((2\gg+1-a)\sqrt{\myfrac{a}{\gg}}\cot\left(\gt\sqrt{a\gg}\right)-m).
\EA$$
If we take $1<a_1<2\gg+1$, we choose $\gt:=\gt_1\in (0,\gt_0]$ such that 
\bel{E15-8}\BA {ll}
(2\gg+1-a_1)\sqrt{\myfrac{a}{\gg}}\cot\left(\gt\sqrt{a_1\gg}\right)>m,
\EA\ee
which implies that $w_{\gt}$ is a supersolution in $S_{\gt}$. In assuming now $a:=a_2>2\gg+1$, we also have, 
$$\BA{lll}
\CP_\gd (w_{\gt})\leq  w_{\gt}'((2\gg+1-a)w_{\gt}'-m)+\ge w_{\gt}\\[2mm]
\phantom{\CP_\gd (w_{\gt})}
\leq |w_{\gt}'|\left((2\gg+1-a)|w_{\gt}'|+m\right)+\ge w_{\gt}\\[2mm]
\phantom{\CP_\gd (w_{\gt})}
\leq |w_{\gt}'|\left(m-(a-2\gg-1)\sqrt{\myfrac{a}{\gg}}\cot\left(\gt\sqrt{a\gg}\right)\right)+\ge w_{\gt}.
\EA$$
We choose $a_2=4\gg+2-a_1$, then 
$$m-(a-2\gg-1)\sqrt{\myfrac{a}{\gg}}\cot\left(\gt\sqrt{a\gg}\right)\leq -c<0\qquad\forall \gt\in (0,\gt_1].
$$
Therefore
$$\CP_\gd (w_{\gt})\leq w_\gt\left(\ge+\myfrac{w'_\gt}{w_\gt}\right).
$$
Since $w'_{\gt_1}<0$ and $w_{\gt_1}({\gt_1})=0^+$, there exists $\ge_0>0$, such that for any $\ge\in (0,\ge_0]$
$$\ge+\myfrac{w'_{\gt_1}(\gr)}{w_{\gt_1}(\gr)}\leq -1\qquad\forall\gr\in (0,\gt_1].
$$
Therefore $w_{\gt_1,a_1}$ and $w_{\gt_1,a_2}$ are respectively supersolution and subsolution of $\CP_\gd(u)=0$ in $S_{\gt_1}$. We extend them in $S'_{\gt_1}$ as smooth functions $\tilde w_{\gt_1,a_1}$ and $\tilde w_{\gt_1,a_2}$ in order $\abs{\tilde \CP_\gd(\tilde w_{\gt_1,a_j})}$ to remain bounded by some constant $M$. Finally $\bar w=\tilde w_{\gt_1,a_1}+M\ge^{-1}$ is a supersolution and 
$\underline w=\tilde w_{\gt_1,a_2}+M\ge^{-1}$ is a subsolution of $\CP_\gd(u)=0$.\smallskip

Next, we replace $\underline w$ by 
$$\underline w_h(\gd)=\underline w(\gd+h)
$$
and $\bar w$ by
$$\bar w_h=\bar w(\gd-h)
$$
for $h$ small enough, we still have a sub and a super solution of $\CP_\gd(u)=0$ in $S_{\gt_1}$ and 
$S_{\gt_1}\setminus S_{h}$.
In the remaining part of $S$, we extend smoothly $\underline w_h$ and $\bar w_h$ in order $\tilde \CP_\gd(\underline w_h)$
and $\tilde \CP_\gd(\bar w_h)$ be bounded.  We can adjust $M$ in order $\CP_\gd(\underline w_h)\leq 0$
and $\CP_\gd(\bar w_h)\geq 0$ in whole $S$, and all these manipulations can be done uniformly with respect to 
$h$ and $\ge$. If $w$ is any $C^2$ solution of $(\ref{E7})$, we prove that it dominates the subsolution $\underline w_h$ in $S$: 
actually, if we assume that $\underline w_h$ and $ w$ are not ordered in $S$, there exists $\gs_0\in S$ such that
$$\underline w_h(\gs_0)- w(\gs_0)=\max\{\underline w_h(\gs)- w(\gs):\gs\in S\}>0.
$$
Since the two functions are $C^2$, 
$$\nabla\underline w_h(\gs_0)=\nabla w(\gs_0)\quad\text{and }\;D^2\underline w_h(\gs_0)\leq D^2 w(\gs_0),
$$
where $D^2$ is the Hessian form, in the sense of quadratic forms, i.e.
$$D^2\underline w_h(\gs_0)(\nabla\underline w_h(\gs_0),\nabla\underline w_h(\gs_0))
\leq D^2 w(\gs_0)(\nabla  w(\gs_0),\nabla w(\gs_0)).
$$
This implies $\CP_\gd(\underline w_h)(\gs_0)>\CP_\gd( w)(\gs_0)=0$, contradiction. Therefore
\bel{E18}\BA {ll}
\underline w_h\leq w\qquad\text{in }\, S,
\EA\ee
uniformly with respect to $h$. Similarly
\bel{E19}\BA {ll}
\bar w_h\geq w\qquad\text{in }\, S.
\EA\ee
Letting $h$ tend to $0$ the claim follows. \qeda
%%%%%%%%%%%%%%%%%%%%%%%%%%%%%%%%%%%%%%%%%%%%%%%%%%%%%%%%%%%%%%%%%%%%%%%%%%%%%%%%%%%%%%%%%%%%%%%%%%%%%%%%%%%%%%%%%%%%%%%%%%%%%%%%%%%%%%%%%%%%%%%%%%%%%%%%%%%%%%%%%%%%%%%%%%%%%%%%%%%%%%%%%%%%%%%%%%%%%%%%%%%%%%%%%%%%%%%%%%%%%%%%%%%%%%%%%%%%%%%%%%%%%%%%%%%%%%%%%%%%%%%%%%%%%%%%%%%%%%%%%%%%%%%%%%%%%%%%%%%%%%%%%%%%%%%%%%%%%%%%%%%%%%%%%%%%%%%%%%%%%%%%%%%%%%%%%%%%%%%%%%%%%%%%%%%%%%%%%%%%%%%%%%%%%%%%%%%%%%%%%%%%%%%%%%%%%%%%%%%%%%%%%%
\subsection{Gradient estimates}

If $\gs_0\in S^{N-1}$ and $R<\gp$, we set $B_R(\gs_0)=\{\gs\in S^{N-1}:\ell (\gs,\gs_0)<R\}$.

\bprop{gradest} Let $0\leq \gd,\ge\leq 1$ and $ w$ be a smooth solution of 
\bel{L1}\BA {ll}
-\myfrac{1}{2}\nabla'|\nabla' w|^2.\nabla w-\gd\Gd w+\gg|\nabla' w|^4+(2\gg+1)|\nabla' w|^2+\ge w=0\quad\text{ in }\,B_R(\gs_0)\subset S,
\EA\ee
where $S$ is a domain of $S^{N-1}$. Then there exists $c=c(N)>0$  such that 
\bel{L1-1}\BA {ll}
\abs{\nabla' w(\gs_0)}\leq\myfrac{c}{\gg R}.
\EA\ee
\es
\nind\Proof
We set $z=|\nabla w|^2$, then $2\Gd_\infty w=\nabla'|\nabla' w|^2.\nabla' w=\nabla' z.\nabla' w$. We define the linearized operator of $\Gd_\infty$ at $ w$ following $h$ by
$$B_ w(h):=\myfrac{d}{dt}\Gd_\infty( w+th)\lfloor_{t=0}=\myfrac{1}{2}\nabla'h.\nabla'z+\nabla' w.\nabla'(\nabla' w.\nabla'h).
$$
Thus the linearized operator of $\Gd_\infty+\gd\Gd$ at $ w$ following $h$ is
\bel{L2}\BA {ll}
\CL_ w(h)=B_ w(h)+\gd\Gd h.
\EA\ee
Thus
$$\CL_ w (z)=\myfrac{1}{2}|\nabla' z|^2+\nabla' w.\nabla'(\nabla'.\nabla z)+\gd\Gd z.
$$
We can re-write $(\ref{L1})$ under the form
\bel{L3}\BA {ll}
\nabla' w.\nabla z=2\left(\gg z^2+(2\gg+1)z+\ge w-\gd\Gd w\right).
\EA\ee
Hence
$$\nabla'(\nabla' w.\nabla'z)=2\left((2\gg z+2\gg+1)\nabla'z+\ge\nabla' w-\gd\nabla'(\Gd w)\right),
$$
and then
$$\nabla w.\nabla'(\nabla' w.\nabla'z)=2\left((2\gg z+2\gg+1)\nabla'z.\nabla' w+\ge z-\gd\nabla'(\Gd w).\nabla' w\right).
$$
By the Weitznb\"ock formula, since ${\rm Ricc}\, (S^{N-1})=(N-2) g_0$ ($g_0$ is the metric tensor on $S^{N-1}$), 
we have
$$\BA {lll}
\myfrac{1}{2}\Gd z=\abs{D^2 w}^2+\nabla'(\Gd w).\nabla' w+(N-2)|\nabla' w|^2\\
\phantom{\myfrac{1}{2}\Gd z}=\abs{D^2 w}^2+\nabla'(\Gd w).\nabla' w+(N-2)z.
\EA$$
Hence
$$\BA {lll}\CL_ w(z)=\myfrac{1}{2}\abs{\nabla z}^2+2\left((2\gg z+2\gg+1)\nabla' z.\nabla' w
+2\ge z-2\gd\nabla' w.\nabla'(\Gd w)\right)\\[2mm]\phantom{\CL_ w(z)-----}
+2\gd\abs{D^2 w}^2+2\gd\nabla'(\Gd w).\nabla' w+2\gd (N-2)z.
\EA$$
Expanding the above identity, we see that the terms of order 3 disappear, hence
\bel{L4}\BA {ll}
\CL_ w(z)=\myfrac{1}{2}\abs{\nabla z}^2+2\left((2\gg z+2\gg+1)\nabla' z.\nabla' w
+2\ge z\right)+2\gd\abs{D^2 w}^2+2\gd (N-2)z.
\EA\ee
If $\xi\in C_c^2(\overline B_R(\gs_0))$, we set $Z=\xi^2z$ and we derive
$$\BA {lll}\CL_ w(Z)=B_ w(\xi^2z)+\gd\Gd(\xi^2z)\\
\phantom{\CL_ w(Z)}
=\xi^2\CL_ w(z)+z\CL_ w\xi^2+2(\nabla' w.\nabla'\xi^2)(\nabla' w.\nabla z)+2\gd\nabla'\xi^2.\nabla z,
\EA$$
where
$$\BA {lll}
\CL_ w\xi^2=\myfrac{1}{2}\nabla' z.\nabla'\xi^2+\nabla' w.\nabla'(\nabla' w.\nabla' \xi^2)+\gd\Gd\xi^2
\\[2mm]
\phantom{\CL_ w\xi^2}
=\myfrac{1}{2}\nabla' z.\nabla'\xi^2+D^2 w(\nabla' \xi^2).\nabla'  w+D^2\xi^2(\nabla'  w).\nabla'  w+\gd\Gd\xi^2\\[2mm]
\phantom{\CL_ w\xi^2}
=\nabla' z.\nabla'\xi^2+D^2\xi^2(\nabla'  w).\nabla'  w+\gd\Gd\xi^2\\[2mm]
\phantom{\CL_ w\xi^2}
\geq \nabla' z.\nabla'\xi^2-z\abs {D^2{\xi^2}}+\gd\Gd\xi^2.
\EA$$
By Schwarz inequality,  $(\Gd w)^2\leq \frac{1}{N-1}\abs{D^2 w}^2$, we derive from  $(\ref{L2})$ and $(\ref{L3})$,
$$\BA {lll}
\CL_ w(z)\geq\myfrac{1}{2}\abs{\nabla z}^2
+4\left(2\gg z+2\gg+1\right)\left(\gg z^2+(2\gg+1)z+\ge w-\gd\Gd w\right)+4\ge z\\[2mm]
\phantom{\CL_ w(z)----}
+\myfrac{2\gd}{N-1}(\Gd w)^2+2\gd (N-2)z\\[2mm]
\phantom{\CL_ w(z)}
\geq \myfrac{1}{2}\abs{\nabla z}^2+\myfrac{\gd}{N-1}(\Gd w)^2+4\gg^2z^3-c_0,
\EA$$
for some $c_0=c_0(N,\gg)>0$. In the sequel the different positive constants $c_j$  which will appear bellow depend only on $N$ and $\gg$. This implies
\bel{L5}\BA {ll}
\CL_ w(Z)\geq z\left(\nabla' z.\nabla'\xi^2-z\abs {D^2{\xi^2}}+\gd\Gd\xi^2\right)+2(\nabla' w.\nabla'\xi^2)(\nabla' w.\nabla z)+2\gd\nabla'\xi^2.\nabla z
\\[2mm]\phantom{\CL_ w(Z)-----}
+\xi^2\left(\myfrac{1}{2}\abs{\nabla z}^2+\myfrac{\gd}{N-1}(\Gd w)^2+4\gg^2z^3-c_0\right).
\EA\ee
We choose $\xi$ such that $0\leq \xi\leq 1$, $\abs{\nabla\xi}\leq c_1R^{-1}$ and $\abs{D^2\xi}\leq c_1R^{-2}$, then 
$$(z+2\gd)\abs{\nabla' z.\nabla'\xi^2}\leq c_1\myfrac{(z+2\gd)\xi}{R}\abs{\nabla' z}
\leq \myfrac{\xi^2}{8}\abs{\nabla' z}^2+c_2\myfrac{(z+2\gd)^2}{R^2},
$$
$$\abs{z\left(\gd\Gd\xi^2-z\abs {D^2{\xi^2}}\right)}\leq \myfrac{c_3(z+2\gd)^2}{R^2},
$$
$$\abs{\nabla' w.\nabla z}\leq \sqrt z\abs{\nabla z},
$$
$$\abs{\nabla' w.\nabla \xi^2}\leq 2\xi\abs{\nabla  w}\abs{\nabla \xi}\leq \myfrac{2c_1\xi\sqrt z}{R},
$$
$$\abs{2(\nabla' w.\nabla'\xi^2)(\nabla' w.\nabla z)}\leq \myfrac{4c_1\xi z\abs{\nabla z}}{R}
\leq \myfrac{\xi^2}{8}\abs{\nabla' z}^2+c_4\myfrac{z^2}{R^2}.
$$
We consider a point $z_0\in B_R$ where $Z$ is maximal, then $\CL_ w(Z)(z_0)\leq 0$, which implies that at this point, 
\bel{L6}\BA {ll}
\xi^2\left(\myfrac{1}{2}\abs{\nabla z}^2+\myfrac{\gd}{N-1}(\Gd w)^2+4\gg^2z^3-c_0\right)
\leq \myfrac{\xi^2}{4}\abs{\nabla z}^2+\myfrac{c_5(z+2\gd)^2}{R^2}.
\EA
\ee
We assume $R\leq 1$ and $2\gd\leq 1$, we multiply by $\xi^4$ and obtain
 \bel{L7}\BA {ll}
\myfrac{1}{4}\abs{\xi^3\nabla z}^2+\myfrac{\gd\xi^6}{N-1}(\Gd w)^2+4\gg^2(\xi^2z)^3
\leq\myfrac{c_6((\xi^2z)^2+1)}{R^2}+c_0.
\EA
\ee
From the inequality 
$$4\gg^2(\xi^2z)^3\leq \myfrac{c_6(\xi^2z)^2}{R^2}+\myfrac{c_7}{R^2},
$$
we deduce 
 \bel{L8}\BA {ll}
\xi^2z\leq \myfrac{c_8}{R^2}\qquad\text{with }c_8=\max\left\{c_7,\frac{c_6}{\gg^2}\right\}.
\EA
\ee
If we assume that $\xi(\gs_0)=1$, we finally infer
 \bel{L8+}\BA {ll}
\abs{\nabla w(\gs_0)}\leq \myfrac{\sqrt{c_8}}{R},
\EA
\ee
which is the claim.
\qeda\medskip

As an immediate consequence, we have

\bcor{estigrad2} Let $0\leq \ge,\gd\leq 1$. If $ w$ is a solution of $(\ref{E7})$ in $S$, it satisfies
 \bel{L9}\BA {ll}
\abs{\nabla w(\gs)}\leq \myfrac{c}{\gg\gr(\gs)}\qquad\forall \gs\in S,
\EA
\ee
for some $c>0$ depending only of $N$.
\es

%%%%%%%%%%%%%%%%%%%%%%%%%%%%%%%%%%%%%%%%%%%%%%%%%%%%%%%%%%%%%%%%%%%%%%%%%%%%%%%%%%%%%%%%%%%%%%%%%%%%%%%%%%%%%%%%%%%%%%%%%%%%%%%%%%%%%%%%%%%%%%%%%%%%%%%%%%%%%%%%%%%%%%%%%%%%%%
\subsection{Proof of Theorem A}
We write $w=w_{\gd,\ge,\gg}$ and
\bel{G1-}\BA {ll}
-\gd\Gd w-\myfrac{1}{2}\nabla\abs{\nabla' w}^2\nabla' w+\gg\abs{\nabla' w}^4+(2\gg+1)\abs{\nabla' w}^2+\ge w=0,
\EA\ee
then $w_{\gd,\ge,\gg}$ satisfies the estimate $(\ref{L9})$. By \rprop{supsub} it satisfies also
\bel{G1-0}\BA {ll}
-\myfrac{1}{\gg}\ln\gr-\myfrac{M}{\ge}\leq w_*\leq w_{\gd,\ge,\gg}\leq w^*+\myfrac{M}{\ge}\leq -\myfrac{1}{\gg}\ln\gr+\myfrac{M}{\ge}.
\EA\ee
The set of functions $\{w_{\gd,\ge,\gg}\}_{\ge,\gd}$ is clearly locally equicontinuous in $S$. By classical stability results on viscosity solutions (see e.g. \cite[Chap 3]{Kat}),  there exist a subsequence $\{w_{\gd_n,\ge,\gg}\}$ and a function $w_{\ge,\gg}$ such that $w_{\gd_n,\ge,\gg}\to w_{\ge,\gg}$, and $w_{\ge,\gg}$ is a viscosity solution of  
\bel{G1-1}\BA {ll}
-\myfrac{1}{2}\nabla'|\nabla' w|^2.\nabla' w+\gb|\nabla' w|^4+(2\gb+1)|\nabla' w|^2+\ge w=0\quad&\text{in }\; S
\\\phantom{-\myfrac{1}{2}\nabla'|\nabla' w|^2.\nabla' w+\gb|\nabla' w|^4+(2\gb+1)|\nabla' w|^2+\ge }
w=\infty&\text{in }\; \prt S. 
\EA\ee
Furthermore $w_{\ge,\gg}$ satisfies the same estimates $(\ref{L9})$ and $(\ref{G1-0})$ as $w_{\gd,\ge,\gg}$. 
Put $\tilde w_{\ge,\gg}(\gs)= w_{\ge,\gg}(\gs)- w_{\ge,\gg}(\gs_0)$ with $\gs_0\in\Gw$, then $\tilde w:=\tilde w_{\ge,\gg}$ satisfies
\bel{G2}\BA {ll}
-\myfrac{1}{2}\nabla\abs{\nabla'\tilde w}^2\nabla'\tilde w+\gg\abs{\nabla'\tilde w}^4+(2\gg+1)\abs{\nabla'\tilde w}^2+\ge\tilde w+\ge w(\gs_0)=0.
\EA\ee
Moreover
$$\abs{\tilde w_{\ge,\gg}(\gs)}=\abs{ w_{\ge,\gg}(\gs)- w_{\ge,\gg}(\gs_0)}\leq \max\left\{\myfrac{c}{\gg\gr(\gt)}:\gt\in [\gs,\gs_0]\right\}|\gs-\gs_0|.
$$
Thus, as $\ge\to 0$, $\ge\tilde w_{\ge,\gg}\to 0$ locally uniformly in $S$. Up to some subsequence $\{\ge_n\}$, $ \tilde w_{\ge_n,\gg}\to  w_{\gg}$ locally uniformly in $S$ and 
$\ge_n  w_{\ge_n,\gg}(\gs_0)\to \gl(\gg)$. As in \cite{PoVe} the expression $\gl(\gg)$ does not depend on $\gs_0$. By analogy with the semilinear case studied in \cite{LaLi}, this last limit is called {\it the ergodic constant}. Furthermore it is easy to check that there exist positive constants $M_1$ and $M_2$ such that $w_1$ and $w_2$ defined by 
\bel{G2'}\BA {ll}
w_1(x)=-\myfrac{1}{\gg}\ln\gr+M_1\gr+M_2\quad\text{ and }\;w_2(x)=-\myfrac{1}{\gg}\ln\gr-M_1\gr-M_2
\EA\ee
are respectively a supersolution and subsolution of $(\ref{G2})$ in $S$ and that there holds
\bel{G2''}\BA {ll}
w_2(x)\leq\tilde w_{\ge,\gg}(x)\leq w_1(x)\qquad\forall x\in S.
\EA\ee
By the same stability results of viscosity solutions, we infer that $ w_{\gg}$ is a positive solution of 
\bel{G3}\BA {ll}
-\myfrac{1}{2}\nabla\abs{\nabla' w}^2\nabla' w+\gg\abs{\nabla' w}^4+(2\gg+1)\abs{\nabla' w}^2+\gl(\gg)=0\quad\text{in }S\\
\phantom{-\myfrac{1}{2}\nabla\abs{\nabla' }^2\nabla' w+\gg\abs{\nabla' w}^4+(2\gg+1)\abs{\nabla' w}^2+\gl(\gg)}
 w=\infty\quad\text{on }\prt S.
\EA\ee
Furthermore, there holds from $(\ref{G2''})$ and $(\ref{L9})$,
\bel{G4}\BA {ll}
\abs{ w_{\gg}+\myfrac{1}{\gg}\ln\gr}\leq M,
\EA\ee
and 
\bel{G5}\BA {ll}
\abs{\nabla w_{\gg}}\leq \myfrac{c}{\gg\gr}.
\EA\ee
\\

\bprop{ergodlim} For any $C^3$ domain $S\subset S^{N-1}$, the ergodic constant $\gl(\gg):= \gl(\gg,S)$ is uniquely determined by $\gg$. Furthermore it is a continuous decreasing function of $\gamma$ and $S$ for the order relation
of inclusion. 
\es

\nind\Proof Assume that the set $\{\ge w_\ge(\gs_0)\}$ of values of the  solutions of $(\ref{G3})$ at  $\gs_0$ admits two different cluster points $\gl_1$ and $\gl_2$.
Then there exist two locally Liptchitz continuous functions $ w_1$ and $ w_2$ satisfying 
\bel{G5-2}\BA {ll}
-\myfrac{1}{2}\nabla'|\nabla' w_i|^2.\nabla' w_i+\gg|\nabla' w_i|^4+(2\gg+1)|\nabla' w_i|^2+\gl_i=0\qquad\text{in } S,
\EA\ee
in the viscosity sense, and such that 
\bel{G5-3} w_i(\gs)=-\myfrac{1}{\gg}\ln\gr(\gs)\left(1+o(1)\right)\quad\text{ as }\gr(\gs)\to 0.
\ee
We can assume that $\gl_1>\gl_2$. For $\ge>0$ let $v=(1+\ge) w_2$. Then
\bel{G5-4}\BA {ll}
-\myfrac{1}{2}\nabla'|\nabla'v|^2.\nabla'v+(1+\ge)^{-1}\gg|\nabla'v|^4+(1+\ge)(2\gg+1)|\nabla'v|^2+(1+\ge)^{3}\gl_2=0\qquad\text{in } S.
\EA\ee
For $X>0$, we put 
$$f(X)=\myfrac{\gg\ge}{1+\ge} X^2-(2\gg+1)\ge X+\gl_1-(1+\ge)^{3}\gl_2.$$
Then
$$f(X)\geq f(X_0)=f\left(\myfrac{(2\gg+1)(1+\ge)}{2\gg}\right)=-\myfrac{\ge(1+\ge)(2\gg+1)^2}{4\gg}+\gl_1-(1+\ge)^{3}\gl_2.
$$
Therefore there exists $\ge_0>0$ such that for any $X\geq 0$, $f(X)\geq 0$, or equivalently
\bel{G5-5}\BA {ll}
(1+\ge)^{-1}\gg X^2+(1+\ge)(2\gg+1)X+(1+\ge)^{3}\gl_2\leq \gg X^2+(2\gg+1)X+\gl_1.
\EA\ee
This implies that
\bel{G5-6}\BA {ll}
-\myfrac{1}{2}\nabla'|\nabla'v|^2.\nabla'v+\gg|\nabla'v|^4+(2\gg+1)|\nabla'v|^2+\gl_1\geq 0\qquad\text{in } S,
\EA\ee
in the viscosity sense. Since $ w_1< v$ near $\prt S$, it follows from comparison principle that $ w_1< v$ in 
$S$. Letting $\ge\to 0$ yields 
\bel{G5-7}\BA {ll}
 w_1\leq w_2\quad\text{in }S.
\EA\ee
Since for any $k\in\BBR$, $ w_1+k$ satisfies the same equation as $ w_1$ and the same estimate $(\ref{G5-3})$ as the $ w_i$ we obtain a contradiction. Thus $\gl=\gl(\gg)$ is uniquely determined. \smallskip

For proving monotonicity, assume $\gg_1>\gg_2>0$ and let $ w_{\ge,1}$ and $ w_{\ge,2}$ be solutions of 
\bel{G5-8}\BA {ll}
-\myfrac{1}{2}\nabla'|\nabla' w_{\ge,i}|^2.\nabla' w_{\ge,i}+\gg_i|\nabla' w_{\ge,i}|^4+(2\gg_i+1)|\nabla' w_{\ge,i}|^2+\ge w_{\ge,i}=0\qquad\text{in } S,
\EA\ee
such that
\bel{G5-9} w_{\ge,i}(\gs)=-\myfrac{1}{\gg_i}\ln\gr(\gs)\left(1+o(1)\right)\quad\text{ as }\gr(\gs)\to 0.
\ee
Then 
$$
-\myfrac{1}{2}\nabla'|\nabla' w_{\ge,1}|^2.\nabla' w_{\ge,1}+\gg_2|\nabla' w_{\ge,1}|^4+(2\gg_2+1)|\nabla' w_{\ge,1}|^2+\ge w_{\ge,1}\leq 0.
$$
Since $ w_{\ge,1}\leq w_{\ge,2}$ near $\prt S$, it follows by comparison principle that $ w_{\ge,1}\leq w_{\ge,2}$ in $S$ and in particular $\ge w_{\ge,1}\leq\ge w_{\ge,2}$. Since $\gl_1=\lim_{n\to\infty}\ge w_{\ge_n,1}(x_0)$ and $\gl_2=\lim_{n\to\infty}\ge w_{\ge_n,2}(x_0)$, we infer that $\gl_1\leq\gl_2$.  \smallskip

For proving the continuity, let $\{\gg_n\}$ be a sequence converging to $\gg$ and let $ w_{n}$  be corresponding solutions of 
\bel{G5-10}\BA {ll}
-\myfrac{1}{2}\nabla'|\nabla' w_{n}|^2.\nabla' w_{n}+\gg_n|\nabla' w_{n}|^4+(2\gg_n+1)|\nabla' w_{n}|^2+\gl(\gg_n)=0\qquad\text{in } S,
\EA\ee
subject to 
\bel{G5-11}\abs{ w_{n}(\gs)+\myfrac{1}{\gg_n}\ln\gr(\gs)}\leq K,
\ee
for some $K>0$ independent of $n$. Since $\{ w_n\}$ is locally bounded in $W^{1,\infty}_{loc}(\Gw)$ we can extract sequences, denoted by
$\{ w_{n_k}\}$, $\{\gl( w_{n_k})\}$ such that $\gl( w_{n_k})\to \bar \gl$ and $ w_{n_k}$ converges locally uniformly to a viscosity solution $ w$ of 
\bel{G5-12}\BA {ll}
-\myfrac{1}{2}\nabla'|\nabla' w|^2.\nabla' w+\gg|\nabla' w|^4+(2\gg+1)|\nabla' w|^2+\bar\gl=0\qquad\text{in } S,
\EA\ee
subject to $ w(\gs)=-\myfrac{1}{\gg}\ln\gr(\gs)\left(1+o(1)\right)$  as $\gr(\gs)\to 0$. The existence of such a function implies that $\bar\gl=\gl(\gg)$. Thus the whole sequence $\{\gl(\gg_n)\}$ converges to $\gl(\gg)$, a fact which implies the continuity. \smallskip

Next, let $S_1\subset S_2$ be two $C^3$ subdomains of $S^{N-1}$. We denote by $w_{\gd,\ge,\gg,S_j}$, $j=1,2$, the solutions of $(\ref{E7})$ respectively in $S_1$ and $S_2$. Since these solutions are limit of solutions with finite boundary values and that the maximum principle holds, we infer that $w_{\gd,\ge,\gg,S_2}\leq w_{\gd,\ge,\gg,S_1}$ in $S_1$. Letting $\gd\to 0$  yields  $w_{\ge,\gg,S_2}\leq w_{\ge,\gg,S_1}$. Taking $\gs_0\in S_1$ and since the ergodic constant is uniquely determined, we have $\ge w_{\ge,\gg,S_2}(\gs_0)\leq \ge w_{\ge,\gg,S_1}(\gs_0)$ and thus $\gl(\gg,S_2)\leq \gl(\gg,S_1)$.
\qeda
\subsection{Proof of Theorem B}

We prove below the following proposition using the result of Theorem C, which proof does not depend on the previous constructions.

\bprop{existence} For any $C^3$ domain $S\subset S^{N-1}$, there exists a unique $\beta:=\gb_s$ such that $\gl(\beta)=\beta+1$. Furthermore $\gb_s$ is a decreasing function of $S$ for the order relation between spherical domains.
\es
\nind\Proof The function $\gamma\mapsto \gl(\gg,S)-\gg$ is continuous and decreasing. For $\ge>0$ we consider two spherical caps $S_i\subset S\subset S_e$; by \rprop{ergodlim}
\bel{G9}\BA {ll}
\gl(\gg,S_e)\leq \gl(\gg,S)\leq \gl(\gg,S_i),
\EA\ee
then
\bel{G10}\BA {ll}
\gl(\gg,S_e)-\gg\leq \gl(\gg,S)-\gg\leq \gl(\gg,S_i)-\gg.
\EA\ee
By Theorem C, there exists $\gg=\gb_{s_e}$ and $\gg=\gb_{s_i}$ such that 
$$\gl(\gb_{s_e},S_e)-\gb_{s_e}=1\quad\text{and }\,\gl(\gb_{s_i},S_i)-\gb_{s_i}=1,
$$
and $\gb_{s_e}<\gb_{s_i}$ unless  $\gl(\gb_{s_e},S_e)=\gl(\gb_{s_i},S_i)$ and $S_i=S_e$.  This implies that
\bel{G11}\BA {ll}
\gl(\gb_{s_e},S)-\gb_{s_e}\geq 1\quad\text{and }\,\gl(\gb_{s_i},S)-\gb_{s_i}\leq 1.
\EA\ee
By continuity there exists a unique $\gb=\gb_s\in [\gb_{s_e},\gb_{s_i}]$ such that $\gl(\gb_s,S)-\gb_s= 1$. To this exponent $\gb$ corresponds a locally Lipschitz continuous function $w$ solution of problem $(\ref{E6})$. Then $\psi=e^{-\gb w}$ is a viscosity solution of $(\ref{E5})$. Notice also that the construction of $\gb_{s}$ and the monotonicity of $S\mapsto \gl(\gg,S)$ imply that $S\mapsto \gb_{s}$ is decreasing. \smallskip

Similarly we can consider separable infinity harmonic functions under the form $(\ref{I2})$ with negative $\gb<0$. We set $\tilde \gb=-\gb$, then $(\ref{E5})$ is replaced by 
\bel{G12}\BA {ll}
-\myfrac{1}{2}\nabla'|\nabla' \psi|^2.\nabla' \psi=\gm(2\gm-1)|\nabla' \psi|^2w+\gm^3(\gm-1)\psi^3\quad&\text{in }\; S
\\\phantom{-\myfrac{1}{2}\nabla'|\nabla' \psi|^2.\nabla' }
\psi=0&\text{in }\; \prt S. 
\EA\ee
If $\psi$ is a positive solution of $(\ref{G12})$, we set
$$ w=-\frac{1}{\gm}\ln\psi.
$$
Then $w$ satisfies
\bel{G13}\BA {ll}\displaystyle
-\myfrac{1}{2}\nabla'\abs{\nabla'w}^2.\nabla'w+\gm\abs{\nabla'w}^4+(2\gm-1)\abs{\nabla'w}^2
=\gm-1\quad\mbox{ in }S
 \\[2mm]\phantom{-----------,-----}\displaystyle
 \lim_{ \gr(\gs)\to 0}w( \gs)=\infty,
\EA\ee
This equation is treated similarly as $(\ref{E6})$. 
\qeda\medskip

\nind\Remark It is an open problem whether the positive functions which satisfy $(\ref{E5})$ are unique up to the multiplication by a constant. This is a sharp contrast with the spherical p-harmonic problem with $1<p<\infty$ where uniqueness is proved by the strong maximum principle and Hopf boundary lemma, and this uniqueness result has been extended to  Lipschitz domain  in \cite{GkVe} using the characterization of the p-Martin boundary obtained by \cite{LeNy}, and a sharp version of boundary Harnack principle. See also Section 3.3 for related results. 
Notice that uniqueness holds when the solutions are spherically radial (see Section 3.2). 
%%%%%%%%%%%%%%%%%%%%%%%%%%%%%%%%%%%%%%%%%%%%%%%%%%%%%%%%%%%%%%%%%%%%%%%%%%%%%%%%%%%%%%%%%%%%%%%%%%%%%%%%%%%%%%%%%%%%%%%%%%%%%%%%%%%%%%%%%%%%%%%%%%%%%%%%%%%%%%%%%%%%%%%%%%%%%%%%%%%%

\mysection{The general case}
\subsection{Problem on the circle}

We consider here the special case $N=2$  and $S$ is the circle in \eqref{E5}. For $k\in\BBN_*$, we set
\bel{F0}\BA {ll}
\gb_k=\myfrac{k^2}{2k+1}\quad\text{and }\;\gm_k=\myfrac{k^2}{2k-1}.
\EA\ee

\bprop {circle}For any $k\in \BBN^*$ there exists two $\frac{\gp}{k}$-anti-periodic $C^1$ functions $\psi_k$ and $\gw_k$ positive  on $(0,\frac{\gp}{k})$ such that $x\mapsto |x|^{-\gb_k} \psi_k(\frac{x}{|x|})$ is infinity harmonic and singular in $\BBR^2\setminus\{0\}$ and  $x\mapsto |x|^{\gm_k} \gw_k(\frac{x}{|x|})$ is infinity harmonic and regular in $\BBR^2$.
\es
\nind\Proof We write $\nabla' \psi=\psi_\gs{\bf e}^{\perp}$ with $S^1\sim\BBR/2\gp$. Thus $(\ref{E5})$ becomes
\bel{F1}\BA {ll}
-\psi^2_\gs \psi_{\gs\gs}=\gb^3(\gb+1)\psi^3+\gb(2\gb+1)\psi^2_\gs \psi, \qquad \psi(0)=\psi(2\pi).
\EA\ee
For $\psi\neq 0$ the equation in \eqref{F1} can be written as
$$- \myfrac{\psi^2_\gs}{\psi^2}\myfrac{\psi_{\gs\gs}}{\psi}=\gb(2\gb+1)\myfrac{\psi^2_\gs}{\psi^2}+\gb^3(\gb+1).
$$
We set $Y=\myfrac{\psi_\gs}{\psi}$, then
$Y_\gs+Y^2=\myfrac{\psi_{\gs\gs}}{\psi}$ and
\bel{Fm}\BA {ll}
-Y^2Y_\gs=Y^4+\gb(2\gb+1)Y^2+\gb^3(\gb+1)=(Y^2+\gb^2)(Y^2+\gb(\gb+1)).
\EA\ee
We first search for solutions with $\gb$ such that $\gb(\gb+1)\geq 0$, $\gb\neq 0$.
%We note that this equation admits the constant solution 
Standard computation yields
\bel{F2}\BA {ll}
\left(\myfrac{\gb}{Y^2+\gb^2}-\myfrac{\gb+1}{Y^2+\gb(\gb+1)}\right)Y_\gs=1,
\EA\ee
and this equation is not degenerate and equivalent to $(\ref{F1})$ as long as $Y\neq 0$.
%\noindent {\it Case 1. } If $\gb=-1$ (that is $\tilde\gb_1$ in \eqref{F0}) then  $(\ref{F2})$ becomes
\noindent If $\gb=-1$ (that is $\tilde\gb_1$ in \eqref{F0}) then  $(\ref{F2})$ becomes
\bel{F2-1}\BA {ll}
\myfrac{Y_\gs}{Y^2+1}=-1,
\EA\ee
thus
\bel{F2-2}\BA {ll}
\tan^{-1}Y(\gs)=-\gs\Longrightarrow Y(\gs)=-\tan\gs\Longrightarrow \psi(\gs)=\sin\gs.
\EA\ee
This corresponds to  the fact that the coordinate functions are separable and infinity harmonic.\\
We assume now $\gb(\gb+1)>0$, or equivalently either $\gb>0$ or $\gb<-1$. We fix $\psi(0)=0$ and consider an interval on the right of $0$ where $\psi>0$. From the equation $\psi$ is concave, thus $\psi_\gs(0)>0$. Because of concavity and periodicity $\psi$ must change sign. We assume that $\gs=\ga$ is the first critical point of $\psi$ which is a singular point for $(\ref{F1})$ and $(\ref{F2})$. We integrate $(\ref{F2})$ on a small interval $(\ga,\gs)$ and get
$$\tan^{-1}\left(\myfrac{Y(\gs)}{\gb}\right)-\sqrt{\myfrac{\gb+1}{\gb}}\tan^{-1}\left(\myfrac{Y(\gs)}{\sqrt{\gb(\gb+1})}\right)=\gs-\ga.
$$
Expanding $\tan^{-1}(x)$ near $x=0$ we obtain
\bel{F2-3}Y(\gs)=-\gb\sqrt[3]{3\gb+3}\sqrt[3]{\gs-\ga}\;(1+o(1))\Longrightarrow \psi(\gs)=\psi(\ga)-C(\gb)\psi(\ga)(\gs-\ga)^{\frac{4}{3}}\;(1+o(1)),
\ee
with $C(\gb)=\frac{3^{\frac{4}{3}}}{4}\gb\sqrt[3]{\gb+1}$, and define $Y(\gs)$ for $\gs\in(\ga,2\ga)$ by imposing $Y(\gs)=-Y(2\ga-\gs)$ and continue this process in order to construct a $2\ga$-antiperiodic solution belonging to $C^{1,\frac{1}{3}}(\BBR)$. Since 
$$\left[\tan^{-1}\left(\myfrac{Y}{\gb}\right)-\sqrt{\myfrac{\gb+1}{\gb}}\tan^{-1}\left(\myfrac{Y}{\sqrt{\gb(\gb+1})}\right)\right]^{\gs=\ga}_{\gs=0}=\ga,
$$
with $Y(0)=\infty$, $Y(\ga)=0$, the condition for $\gp$-antiperiodicity is therefore
$$\left(\sqrt{\myfrac{\gb+1}{\gb}}-1\right)\myfrac{\gp}{2}=\ga\Longleftrightarrow \gb=\myfrac{\gp^2}{4(\ga^2+\ga\gp)}.
$$
If $\gb>0$, the periodicity condition yields
\bel{F3}\BA {ll}
\sqrt{\myfrac{\gb+1}{\gb}}=1+\myfrac{1}{k}\Longleftrightarrow\gb=\gb_k=\myfrac{k^2}{1+2k}.
\EA\ee
If $\gb<-1$
$$\left(1-\sqrt{\myfrac{\gb+1}{\gb}}\,\right)\gp=\ga,
$$
and the periodicity condition implies
\bel{F3+}\BA {ll}
\sqrt{\myfrac{\gb+1}{\gb}}=1-\myfrac{1}{k}\Longleftrightarrow\gb=\tilde\gb_k=\myfrac{k^2}{1-2k}.
\EA\ee\smallskip

\nind The case $\gb(\gb+1)<0$, or equivalently $-1<\gb<0$, is easily ruled out. We find that \eqref{Fm} has the constant solution $Y(\gs)=\sqrt{-\gb(\gb+1)}$, meaning
$\psi(\gs)=C\exp(\gs\sqrt{-\gb(\gb+1)})$, which is by no mean periodic.
%%%%%%%%%%%%%%%%%%%%%%%%%%%%%%%%%%%%%%%%%%%%%%%%%%%%%%%%%%%%%%%%%%%%%%%%%%%%%%%%%%%%%%%%%%%%%%%%%%%%%%%%%%%%%%%%%%%%%%%%%%%%%%%%%%%%%%%%%%%%%%%%%%%%%%%%%%%%%%%%%%%%%%%%%%%%%%
On the other hand, in this case we can write $(\ref{F2})$ under the form
\bel{F4}\BA {ll}	
\myfrac{d}{d\gs}\left(\tan^{-1}\left(\myfrac{Y}{\gb}\right)-\myfrac{1}{2}\sqrt{\myfrac{\gb+1}{-\gb}}\ln
\left(\myfrac{|Y-\sqrt{-\gb(\gb+1)}|}{Y+\sqrt{-\gb(\gb+1)}}\right)\right)=1.
\EA\ee
Since $Y$ runs from $Y(0)=\infty$ to $Y(\ga)=0$, there must be a value $\gs_0$ where $Y(\gs_0)=\sqrt{-\gb(\gb+1)}$. We can integrate \eqref{F4} on $(0,\gs_0-\ge)$ and let $\ge\to 0$. Since $\gb<0$, it yields
$$\BA {lllll}
\myfrac{\gp}{2}+\tan^{-1}\left(\myfrac{Y(\gs_0-\ge)}{\gb}\right)
-\myfrac{1}{2}\sqrt{\myfrac{\gb+1}{-\gb}}\ln
\left(\myfrac{|Y(\gs_0-\ge)-\sqrt{-\gb(\gb+1)}|}{Y(\gs_0-\ge)+\sqrt{-\gb(\gb+1)}}\right)=\gs_0-\ge.
\EA$$
The left-hand side expression tends to $\infty$ when $\ge\to 0$, a contradiction. Hence there are no solutions with $\gb\in (-1,0)$. This ends the proof of the proposition.
\qeda\medskip

\nind\Remark When $k=1$ the coordinate functions are infinity harmonic and vanish on a straight line. When $k=2$, the regular solution with $\gm_1=\frac{4}{3}$ is 
$$u(x,y)=x^{\frac{4}{3}}-y^{\frac{4}{3}}.
$$
Its existence is due to Aronsson \cite{Aro1}. The corresponding circular function, $\gw(\gs)=(\cos\gs)^{\frac{4}{3}}-(\sin\gs)^{\frac{4}{3}}$, admits four nodal sets on $S^1$. When $k=1$, then $\gb_1=\frac{1}{3}$. It is proved in \cite{Bat} that any positive infinity harmonic function in a half-space which vanishes on the boundary except at one point blows-up like the separable infinity harmonic function $u(r,\gs)=r^{-\frac{1}{3}}\psi(\gs)$.

\subsection{The spherical cap problem}

\nind{\it Proof of Theorem C}. The following representation of $S^{N-1}$ is classical
$$S^{N-1}=\left\{\gs=(\sin\gf\,\gs',\cos\gf):\gs'\in S^{N-2},\gf\in [0,\gp]\right\}.
$$
Then $\nabla'\psi=\psi_\gf{\bf e}+\nabla'_{\gs'}\psi$ where ${\bf e}$ is a tangent unit downward vector to $S^{N-1}$ following the great circle going through the point $\gs$. Then $|\nabla'\psi|^2=\psi^2_\gf+|\nabla'_{\gs'}\psi|^2$, thus, if $\psi$ depends only on $\phi$, we have
$$\frac12\nabla'|\nabla'\psi|^2.\nabla'\psi=\psi_{\phi\phi}\psi^2_\gf.
$$
Therefore such a function $\psi$, if it is a $C^1$ solution of $(\ref{E5})$ in the spherical cap $S_\ga$ defined for $\gf\in (0,\ga)$, satisfies
\bel{X1}\BA{lll}
-\psi_{\phi\phi}\psi^2_\gf=\gb(2\gb+1)\psi^2_\gf\psi_\gf+\gb^3(\gb+1)\psi\quad\text{in }(0,\ga)\\
\phantom{-'}\,\psi_\gf(0)=0\,,\;\psi(\ga)=0.
\EA\ee
The conclusion follows from \rprop{circle}.\qeda\medskip

\nind\Remark If $\ga=\gp$, the exponent $\gb_+$ is $\frac{1}{8}$ and $\psi:=\psi_{_{_{\Gs^c}}}$ is a positive solution of 
\bel{X2}\BA {ll}
-\psi^2_\gs \psi_{\gs\gs}=\myfrac{9}{4096}\psi^3+\myfrac{5}{32}\psi^2_\gs \psi\quad\text{in }(-\gp,\gp)\\[2mm]
\phantom{\psi_{}}
\psi(-\gp)=\psi(\gp)=0.
\EA\ee
Then the function $u_{_{_{\Gs^c}}}(r,\gs)=r^{-\frac{1}{8}}\psi_{_{_{\Gs^c}}}(\gs)$ is an infinity harmonic function in $\BBR^N\setminus L_{_{_{\Gs^c}}}$, which vanishes on  the half line
$L_{_{_{\Gs^c}}}:=\{x=t\Gs:t\geq 0\}$. The function $Y=\frac{\psi_\gs}{\psi}$ can be computed implicitly on $(0,\gp)$ thanks to the identity
\bel{Y3}\BA {ll}
\tan^{-1}(8Y(\gs))-3\tan^{-1}(\frac 83Y(\gs))=\gs.
\EA\ee
This yields, with $Z=\frac{8Y}{3}$, $\tan^{-1}(3Z(\gs))-3\tan^{-1}(Z(\gs))=\gs$, hence
\bel{Y4}\BA {ll}
\myfrac{3Z-\tan(3\tan^{-1}(Z)))}{1+3Z\tan(3\tan^{-1}(Z))}=\tan\gs,
\EA\ee
since 
$$\tan(3x)=\myfrac{3\tan x-\tan^3x}{1-3\tan^2x}.
$$
This yields
\bel{Y4x}\BA {ll}
\myfrac{-8Z^3}{1+6Z^2-3Z^4}=\tan\gs,
\EA\ee
which gives the value of $Y$ by solving a fourth degree equation and then $\psi=\psi_{_{_{\Gs^c}}}$ by integrating $Y$.\medskip

Using Theorem C we can prove the existence of a singular infinity harmonic function in a cone $C_{_{S_{\gk,\ga}}}$ generated by a spherical annulus $S_{\gk,\ga}$ of the spherical points  with azimuthal angle $\gk<\gf<\ga$. 

\bprop{annulus} Assume $0\leq\gk<\ga<\gp$ and let $\gn=\frac{1}{2}(\ga-\gk)$. Then there exists a positive singular  infinity harmonic function $u_{_{S_{\gk,\ga}},+}$ and a regular infinity harmonic function $u_{_{s_{\gk,\ga}},-}$ in $C_{_{S_{\gk,\ga}}}$ which vanish respectively on $\prt C_{_{S_{\gk,\ga}}}\setminus\{0\}$ and $\prt C_{_{S_{\gk,\ga}}}$ under the form $u_{_{S_{\gk,\ga}},+}(r,\gs)=r^{-\gb_{{s_{\gk,\ga}}}}\psi_{_{S_{\gk,\ga}}}(\gs)$ and $u_{_{S_{\gk,\ga}},-}(r,\gs)=r^{\gm_{_{S_{\gk,\ga}}}}\gw_{_{S_{\gk,\ga}}}(\gs)$ where 
\bel{X3}\BA {ll}\gb_{_{S_{\gk,\ga}}}=\myfrac{\gp^2}{4\gn(\gp+\gn)}\;\text{ and}\quad\gw_{_{S_{\gk,\ga}}}=\myfrac{\gp^2}{4\gn(\gp-\gn)},
\EA\ee
and $\psi_{_{S_{\gk,\ga}}}$ and $\gw_{_{S_{\gk,\ga}}}$  are positive solutions of $(\ref{F1})$ in $S_{\gk,\ga}$  vanishing at $\gk$ and $\ga$ with $\gb=\gb_{_{S_{\gk,\ga}}}$ and  $\gb=-\gm_{_{S_{\gk,\ga}}}$ respectively.  
\es
\nind\Proof By Theorem C there exists a positive and even solution $\tilde\psi$  of 
\bel{X4}\BA {ll}
\tilde\psi_{\phi\phi}\psi^2_\gf=\gb(2\gb+1)\tilde\psi^2_\gf\tilde\psi_\gf+\gb^3(\gb+1)\tilde\psi\quad\text{in }(-\gn,\gn):=\left(\frac{1}{2}(\gk-\ga),\frac{1}{2}(\ga-\gk)\right)\\
\psi(-\gn)=0\,,\;\psi(\gn)=0,
\EA\ee
with $\gb=\gb_{_{S_{\gk,\ga}}}$ or $\gb=-\gm_{_{S_{\gk,\ga}}}$. Then $\phi\mapsto \psi(\gf):=\tilde\psi(\gf+\frac{1}{2}(\gk+\ga))$ is a positive solution of $(\ref{F1})$ in $(\gk,\ga)$. The proof follows.\qeda 
%%%%%%%%%%%%%%%%%%%%%%%%%%%%%%%%%%%%%%%%%%%%%%%%%%%%%%%%%%%%%%%%%%%%%%%%%%%%%%%%%%%%%%%%%%%%%%%%%%%%%%%%%%%%%%%%%%%%%%%%%%%%%%%%%%%%%%%%%%%%%%%%%%%%%%%%%%%%%%%%%%%%%%%%%%%%%%%%%%%%

\medskip
%%%%%%%%%%%%%%%%%%%%%%%%%%%%%%%%%%%%%%%%%%%%%%%%%%%%%%%%%%%%%%%%%%%%%%%%%%%%%%%%%%%%%%%%%%%%%%%%%%%%%%%%%%%%%%%%%%%%%%%%%%%%%%%%%%%%%%%%%%%%%%%%%%%%%%%%%%%%%%%%%%%%%%%%%%%%%%%%%%%%

The next technical lemma is a variant of Theorem C and \rprop{annulus}.

\blemma{punct} Assume $0<\ga<\gp$ and $\ge,\gg>0$. Then the solution $v=v_{\ge,\gg,\ga}$ of 
\bel{X5-}\BA {ll}
-v'^2v''+\gg v'^4+(2\gg+1)v'^2+\ge v=0\qquad\text{in } (0,\ga)\\\phantom{------,,-}
v(0)=\infty\,,\;v(\ga)=\infty,
\EA\ee
is an increasing function of $\ge$. If $0<\gs_0<\ga$, there exists $\gl=\gl(\ga,\gg)=\lim_{\ge\to 0}\ge v_{\ge,\gg,\ga}(\gs_0)$
and this value is independent of $\gs_0$. The function $\tilde v=\tilde v_{\ge,\gg,\ga}=v_{\ge,\gg,\ga}-v_{\ge,\gg,\ga}(\gs_0)$ converges locally uniformly in $(0,\ga)$ to a solution $v=v_{\gg,\ga}$ of 
\bel{X5}\BA {ll}
-v'^2v''+\gg v'^4+(2\gg+1)v'^2+\gl=0\qquad\text{in } (0,\ga)\\\phantom{-------}
v(0)=\infty\,,\;v(\ga)=\infty,
\EA\ee
with
\bel{X6}\BA {ll}
\gl(\ga,\gg)=\myfrac{1}{4\gg^3}\left(\myfrac{\gp^2}{\ga^2}-\gg(2\gg+1)\right)^2.
\EA\ee
Furthermore
\bel{X6+}\BA {ll}
v_{\gg,\ga}(\gf)=-\myfrac{1}{\gg}\ln\gf\quad\text{as }\;\gf\to 0,
\EA\ee
and
\bel{X6*}\BA {ll}
v_{\gg,\ga}(\gf)=-\myfrac{1}{\gg}\ln(\ga-\gf)\quad\text{as }\;\gf\to \ga.
\EA\ee
\es
\nind
\Proof From the proof of Theorem A, we know that $v_{\ge,\gg,\ga}$ is an increasing function of $\ge$. It satisfies
 estimates $(\ref{G1-0})$ and $(\ref{L9})$. Furthermore there exists $\gl=\gl(\ga,\gg)=\lim_{\ge\to 0}\ge v_{\ge,\gg,\ga}(\gs_0)\geq 0$, 
which is a value independent of $\gs_0\in (0,\ga)$, and $\tilde v=\tilde v_{\ge,\gg,\ga}=v_{\ge,\gg,\ga}-v_{\ge,\gg,\ga}(\gs_0)$   converges locally uniformly in $(0,\ga)$ to a solution $v=v_{\gg,\ga}$ of $(\ref{X5})$. We set $Y=-\gg v'$, then 
\bel{X6++}\BA {ll}
Y^2Y'+Y^4+\gg(2\gg+1)Y^2+\gl\gg^3=0\qquad\text{in } (0,\ga)\\\phantom{-------}
Y(0)=\infty\,,\;Y(\ga)=-\infty.
\EA\ee
We write it under the separable form
$$\left(\myfrac{Y^2}{Y^4+\gg(2\gg+1)Y^2+\gl\gg^3}\right)Y'=-1\Leftrightarrow \left(\myfrac{A^2}{A^2-B^2}\myfrac{1}{Y^2+A^2}-\myfrac{B^2}{A^2-B^2}\myfrac{1}{Y^2+B^2}\right)Y'=-1,
$$
for some $A,B>0$ and with $A>B$ if we assume $2\gg+1>2\gg\gl$. Actually $A^2B^2=\gl\gg^3$ and $A^2+B^2=\gg(2\gg+1)$. Thus
\bel{X7-}\left(A\tan^{-1}\left(\myfrac{Y}{A}\right)-B\tan^{-1}\left(\myfrac{Y}{B}\right)\right)'=B^2-A^2.
\ee
By integration on $(0,\ga)$ we derive the identity
\bel{X7}
A+B=\myfrac{\gp}{\ga}.
\ee
Since $A+B=\sqrt{\gg(2\gg+1)+2\sqrt{\gl\gg^3}}$, we deduce   $(\ref{X6})$
from $(\ref{X7})$. Finally, since 
$$\tan^{-1}z=\frac{\gp}{2}-\frac{1}{z}+\frac{1}{3z^3}+O(z^{-5})\quad\text{when }\,z\to\infty,
$$
we derive
$$-\myfrac{1}{\gg v'(\gf)}=\myfrac{1}{Y(\gf)}=\gf +O(\gf^3)\quad\text{when }\,\gf\to 0
$$
from $(\ref{X7-})$, which implies $(\ref{X6+})$ by l'Hospital rule. Relation $(\ref{X6*})$ is proved similarly. 

\qeda

\medskip

Next we denote by $S_\ga(a)$ the spherical cap with vertex $a\in S^{N-1}$ and azimuthal opening $\ga$ from $a$ and $S^*_\ga(a)=S_\ga(a)\setminus\{a\}$. The next statement is a rephrasing of \rlemma{punct} in a geometric framework.

\bcor{punct-cap} Let $\ga$, $\ge$ and $\gg>0$ be as in \rlemma{punct} and $a\in S^{N-1}$. Then there exists a unique solution $w=w_{a,\ga,\gg,\ge}$ of 
\bel{X8-0}\BA {ll}
-\myfrac{1}{2}\nabla'\abs{\nabla' w}^2.\nabla' w+\gg \abs{\nabla' w}^4+(2\gg+1)\abs{\nabla' w}^2+\ge w=0\quad\text{in }S^*_\ga(a)
\\\phantom{-----------------,}\displaystyle
\lim_{\ell (\gs,a)\to 0}w(\gs)=\infty
\\\phantom{-----------------,}\displaystyle
\lim_{\ell (\gs,a)\to \ga}w(\gs)=\infty,
\EA\ee
rotationally invariant with respect to $a$. If $a$ is replaced by $a'\in S^{N-1}$, the solution $w_{a',\ga,\gg,\ge}$ of 
$(\ref{X8-0})$ in $S^*_\ga(a')$ is derived from $w_{a,\ga,\gg,\ge}$ by an orthogonal transformation exchanging $a$ and $a'$. The mapping $\ge\mapsto w_{a,\ga,\gg,\ge}$ is decreasing and for any $\gs_0\in S^*_\ga(a)$
\bel{X9}
\lim_{\ge\to 0}\ge w_{a,\ga,\gg,\ge}(\gs_0)=\gl(\gg,S^*_\ga(a)):=\gl(\gg,S^*_\ga),
\ee
(this notation is coherent with $\gl(\gg,S)$ already used, furthermore its value does not depend on $a$). The function $\tilde w_{a,\ga,\gg,\ge}=w_{a,\ga,\gg,\ge}- w_{a,\ga,\gg,\ge}(\gs_0)$ converges locally uniformly in $S^*_\ga(a')$
 to  the unique  viscosity solution $w:=w_{a,\ga,\gg}$ rotationally invariant with respect to $a$ and vanishing at $\gs_0$ of 
 \bel{X10-0}\BA {ll}
-\myfrac{1}{2}\nabla'\abs{\nabla' w}^2.\nabla' w+\gg \abs{\nabla' w}^4+(2\gg+1)\abs{\nabla' w}^2+\gl(\gg,S^*_\ga)=0\quad\text{in }S^*_\ga(a)
\\\phantom{-------------------,}\displaystyle
\lim_{\ell (\gs,a)\to 0}w(\gs)=\infty
\\\phantom{-------------------,}\displaystyle
\lim_{\ell (\gs,a)\to \ga}w(\gs)=\infty.
\EA\ee
Finally
 \bel{X11}\BA {ll}
 w_{a,\ga,\gg}(\gs)=-\myfrac{1}{\gg}\ln\left(\ell(a,\gs)\right)+O(1)\quad\text{as }\;\gs\to a,
 \EA\ee
 and
  \bel{X12}\BA {ll}
  w_{a,\ga,\gg}(\gs)=-\myfrac{1}{\gg}\ln\left(\ga-\ell(a,\gs)\right)+O(1)\quad\text{as }\;\ell(\gs, a)\to 0.
   \EA\ee
\es

The following statement is formally similar to \rcor{punct-cap}. It makes more precise the approximations used in the proof of Theorem A, in the construction of the proof of Theorem B in the case $\ga=\gp$.

\bcor{punct-sph} Let $\ge$ and $\gg>0$ and $a\in S^{N-1}$. Then there exists a unique rotationally invariant with respect to $a$ solution $v=v_{a,\gg,\ge}$, of 
\bel{X8}\BA {ll}
-\myfrac{1}{2}\nabla'\abs{\nabla' v}^2.\nabla' v+\gg \abs{\nabla' v}^4+(2\gg+1)\abs{\nabla' v}^2+\ge v=0\quad\text{in }S^{N-1}\setminus\{a\}
\\\phantom{-----------------,}\displaystyle
\lim_{\ell (\gs,a)\to 0}v(\gs)=\infty.
\EA\ee
Furthermore, for any $\gs_0\in S^{N-1}\setminus\{a\}$, 
\bel{X9-0}
\lim_{\ge\to 0}\ge v_{a,\gg,\ge}(\gs_0)=\Gl(\gg):=\myfrac{1}{4\gg^3}\left(\myfrac{1}{4}-\gg(2\gg+1)\right)^2.
\ee
The function $\tilde v_{a,\gg,\ge}=v_{a,\gg,\ge}- v_{a,\gg,\ge}(\gs_0)$ converges locally uniformly in $S^{N-1}\setminus\{a\}$
 to  the unique  viscosity solution $v:=v_{a,\gg}$ rotationally invariant with respect to $a$ and vanishing at $\gs_0$ of 
 \bel{X10}\BA {ll}
-\myfrac{1}{2}\nabla'\abs{\nabla' v}^2.\nabla' v+\gg \abs{\nabla' v}^4+(2\gg+1)\abs{\nabla' v}^2+\Gl(\gg)=0\quad\text{in }S^{N-1}\setminus\{a\}
\\\phantom{------------------}\displaystyle
\lim_{\ell (\gs,a)\to 0}v(\gs)=\infty.
\EA\ee
Finally
 \bel{X11-1}\BA {ll}
 v_{a,\gg}(\gs)=-\myfrac{1}{\gg}\ln\left(\ell(a,\gs)\right)+O(1)\quad\text{as }\;\gs\to a.
 \EA\ee
\es

As in \rcor{punct-cap}, if $a$ is replaced by $a'\in S^{N-1}$, the solution $v_{a',\gg,\ge}$ of $(\ref{X8})$ in $S^{N-1}\setminus\{a'\}$ is derived from $v_{a,\gg,\ge}$ by an orthogonal transformation exchanging $a$ and $a'$. The mapping $\ge\mapsto v_{a,\gg,\ge}$ is decreasing.\medskip

%%%%%%%%%%%%%%%%%%%%%%%%%%%%%%%%%%%%%%%%%%%%%%%%%%%%%%%%%%%%%%%%%%%%%%%%%%%%%%%%%%%%%%%%%%%%%%%%%%%%%%%%%%%%%%%%%%%%%%%%%%%%%%%%%%%%%%%%%%%%%%%%%%%%%%%%%%%%%%%%%%%%%%%%%%%%%%%%%%%%%%%%%%%%%%%%%%%%%%%%%%%%%%%%%%%%%%%%%%%%%%%%%%%%%%%%%%%%%%%%%%%%%%%%%%%%%%%%%%%%%%%%%%%%%%%%%%%%%%%%%%%%%%%%%%%%%%%

\subsection{Proof of Theorem D}
 {\it Step 1: Approximate solutions}. 
We consider an increasing sequence of smooth spherical domains, $\{S_k\}$ such that 
$$\displaystyle S_k\subset\overline S_k\subset S_{k+1}\subset S\,\text{ and }\;\bigcup_k S_k=S,$$
 To each domain we associate  the positive exponent $\gb_k:=\gb_{s_k}$ and the corresponding spherical infinity-harmonic function $\psi_k:=\psi_{s_k}$ defined in $S_k$ and such that $\psi_k(\gs_0)=1$ for some $\gs_0\in S_1$, so that the function 
$u_{k}(r,.)=r^{-\gb_{k}}\psi_{k}$ is infinity-harmonic in 
the cone $C_{S_k}$ and vanishes on $\prt C_{S_k}\setminus\{0\}$. For $\gg,\gd,\ge>0$, we denote by $w_{k,\gg,\gd,\ge}$  the solution of
\bel{Y5}\BA {ll}
-\gd\Gd w-\myfrac{1}{2}\nabla'\abs{\nabla'w}^2.\nabla'w+\gg\abs{\nabla'w}^4+(2\gg+1)\abs{\nabla'w}^2+\ge w=0\quad\text{in } S_k\\[2mm]\phantom{\gd\Gd w-\myfrac{1}{2}\nabla'\abs{\nabla'w}^2+\gg\abs{\nabla'w}^4+(2\gg+1)\abs{\nabla'w}^2}
\displaystyle\lim_{\gr_k(\gs)\to 0}w(\gs)=\infty,
\EA\ee
where $\gr_k(.)=\dist(.,\prt S_k)$. By the maximum principle the functions
  $w_{k,\gg,\gd,\ge}$ is positive and  the following comparison relations hold:
\bel{Y6+}\BA {llll}
&(i)\qquad\qquad\qquad w_{\ell,\gg,\gd,\ge}\leq w_{k,\gg,\gd,\ge}\quad &\text{ in }S_k\;&\forall k\leq \ell,\qquad\qquad\qquad
\\[4mm]
&(ii)\qquad\qquad\qquad w_{k,\gg,\gd,\ge}\leq w_{k,\gg,\gd,\ge'}\quad &\text{ in }S_k &\forall \ge'\leq \ge,\qquad
\\[4mm]
&(iii)\qquad\qquad\qquad w_{k,\gg,\gd,\ge}\leq w_{k,\gg',\gd,\ge}\quad &\text{ in }S_k &\forall \gg'\leq \gg.\qquad
\EA\ee
Furthermore it follows from \rcor{estigrad2},
\bel{Y7}\BA {ll}
\qquad\qquad \abs{\nabla w_{k,\gg,\gd,\ge}(\gs)}\leq \myfrac{c}{\gr_k(\gs)}\qquad\forall \gs\in S_k,\EA\ee
where $c=c(N)$. Moreover, similarly as in $(\ref{G1-0})$,
\bel{Y8}\BA {llll}
\qquad\;\; &-\myfrac{1}{\gg}\ln\gr_k(\gs)-\myfrac{M_k}{\ge}\leq  w_{k,\gg,\gd,\ge}(\gs)\leq -\myfrac{1}{\gg}\ln\gr_k(\gs)+\myfrac{M_k}{\ge}\qquad\;\forall \gs\in S_k.
\EA\ee
We let $\gd\to 0$ and derive that, up to a subsequence, $w_{k,\gg,\gd_n,\ge}\to w_{k,\gg,\ge}$ locally uniformly in $S_k$. The function $w_{k,\gg,\ge}$ satisfy $(\ref{Y6+})$, $(\ref{Y7})$ and $(\ref{Y8})$ and is a viscosity solution of 
\bel{Y9}\BA {ll}
-\myfrac{1}{2}\nabla'\abs{\nabla'w}^2.\nabla'w+\gg\abs{\nabla'w}^4+(2\gg+1)\abs{\nabla'w}^2+\ge w=0\quad\text{in } S_k\\[2mm]\phantom{\myfrac{1}{2}\nabla'\abs{\nabla'w}^2+\gg\abs{\nabla'w}^4+(2\gg+1)\abs{\nabla'w}^2}
\displaystyle\lim_{\gr_k(\gs)\to 0}w(\gs)=\infty.
\EA\ee
Furthermore the mapping $(k,\ge)\mapsto w_{k,\gg,\ge}$ is nonincreasing, and if we let $k\to\infty$, then $w_{k,\gg,\ge}\downarrow w_{\gg,\ge}$. The function $w_{\gg,\ge}$ is defined in $S$ and is nonincreasing functions of $\ge$ and $\gg$. Furthermore there holds
\bel{Y14}\BA {ll}
(i)\qquad\qquad\qquad  w_{\gg,\ge}\leq w_{k,\gg,\ge}\leq w_{1,\gg,\ge}\qquad &\forall k,\ell\geq 1,
\\[2mm]
(ii)\qquad\qquad\qquad \abs{\nabla w_{\gg,\ge}(\gs)}\leq \myfrac{c}{\gr(\gs)}\qquad&\forall \gs\in S,
\EA\ee
where $c=c(N)$. Estimate $(\ref{Y14})$-(i) can be made more precise in the following way: {\it for each $\gs\in S$, there is $k_\gs\in\BBN$ such that $\gs\in S_{k_\gs}$ and}
\bel{Y14+0}\BA {llll}
\quad\;\; -\myfrac{1}{\gg}\ln\gr(\gs)-\myfrac{1}{\ge}\displaystyle\max_{k\geq k_\gs} M_k\leq  w_{\gg,\ge}(\gs)\leq -\myfrac{1}{\gg}\ln\gr(\gs)+
\myfrac{1}{\ge}\displaystyle\min_{k\geq k_\gs} M_k.
\EA\ee
\medskip
%%%%%%%%%%%%%%%%%%%%%%%%%%%%%%%%%%%%%%%%%%%%%%%%%%%%%%%%%%%%%%%%%%%%%%%%%%%%%%%%%%%%%%%%%%%%%%%%%%%%%%%%%%%%%%%%%%%%%%%%%%%%%%%%%%%%%%%%%%%%%%%%%%%%%%%%%%%%%%%%%%%%%%%%%%%%%%%%%%%%%%%%%%%%%%%%%%%%%%%%%%%%%%%%%%%%

\nind{\it Step 2: Boundary blow-up.} 
%%%%%
The compactness of approximate solutions vanishing at a fixed point in the local uniform convergence topology is easy to obtain thanks to the uniform estimate of the gradient. The main difficulty is to preserve the boundary blow-up when the parameters $k,\gg,\gd,\ge$ tend to their respective limit. \smallskip

\nind{\it Case 1}. We first assume that there exist $\gs_0\in S$, two decreasing sequences $\{\ge_n\}$, $\{\gd_\ell\}$ converging to $0$ and an increasing sequence $\{k_j\}$ tending to infinity with the property that
 \bel{Z1}\BA {ll}
\ge_nw_{k_j,\gg,\gd_\ell,\ge_n}(\gs_0)\leq \ge_m w_{k_j,\gg,\gd_\ell,\ge_m}(\gs_0)\quad\text{for all }\,
m<n\,,\;j,\ell\in\BBN.
\EA\ee
Since $\tilde w_{k_j,\gg,\gd_\ell,\ge_n}=w_{k_j,\gg,\gd_\ell,\ge_n}-w_{k_j,\gg,\gd_\ell,\ge_n}(\gs_0)$ satisfies
 \bel{Z1-1}\BA {ll}
-\gd_\ell\Gd'\tilde w-\myfrac{1}{2}\nabla'\abs{\nabla'\tilde w}^2.\nabla'\tilde w+\gg\abs{\nabla'\tilde w}^4+(2\gg+1)\abs{\nabla'\tilde w}^2
+\ge_n\tilde w+\ge_n w_{k_j,\gg,\gd_\ell,\ge_n}(\gs_0)=0\quad\text{in } S_{k_j}\\
\phantom{-\gd'_\ell\Gd'\tilde w'--------\nabla'\tilde w',',,---_{k'\gg'}----------}
\displaystyle\lim_{\gr_{k_j}(\gs)\to 0}\tilde w(\gs)=\infty,
\EA\ee
there holds
 \bel{Z1-2}\BA {ll}
\tilde w_{k_j,\gg,\gd_\ell,\ge_n}\geq \tilde w_{k_j,\gg,\gd_\ell,\ge_m}\quad\text{for all }\,
m<n\,,\;j,\ell\in\BBN.
\EA\ee
Letting $\gd_\ell\to 0$ we derive that $w_{k_j,\gg,\gd_\ell,\ge_n}\to w_{k_j,\gg,\ge_n}$ locally uniformly in $S_{k_j}$ and $w_{k_j,\gg,\ge_n}$ satisfies
 \bel{Z2}\BA {ll}
-\myfrac{1}{2}\nabla'\abs{\nabla'w}^2.\nabla'w+\gg\abs{\nabla'w}^4+(2\gg+1)\abs{\nabla'w}^2
+\ge_n w=0\quad\text{in }\; S_{k_j}\\
\phantom{----------\nabla'W,,-----}
\displaystyle\lim_{\gr_{k_j}(\gs)\to 0}w(\gs)=\infty.
\EA\ee
Furthermore, for all 
$m<n\,,\;j\in\BBN$,
 \bel{Z1'}\BA {lll}
 (i)\qquad\qquad &w_{k_j,\gg,\ge_n}\geq w_{k_j,\gg,\ge_m}\\[2mm]
 (ii)\qquad\qquad &w_{k_j,\gg,\ge_n}-w_{k_j,\gg,\ge_n}(\gs_0)\geq w_{k_j,\gg,\ge_m}-w_{k_j,\gg,\ge_m}(\gs_0)
 \\[2mm]
(iii)\qquad\qquad&\ge_nw_{k_j,\gg,\ge_n}(\gs_0)\leq \ge_m w_{k_j,\gg,\ge_m}(\gs_0).
\EA\ee
By monotonicity with respect to $S_{k_j}$, $w_{k_j,\gg,\ge_n}\downarrow w_{\gg,\ge_n}$ as $j\to\infty$. Let $a\in\prt S$ and  $v_{a,\gg,\ge_n}$ be the solution of $(\ref{X8})$ with $\ge=\ge_n$ which exists by \rcor{punct-sph}. Then 
 \bel{Z1''}
v_{a,\gg,\ge_n} \leq w_{k_j,\gg,\ge_n} \qquad\text{in }S_{k_j},
 \ee
which yields
 \bel{Z1'''}
v_{a,\gg,\ge_n} \leq w_{\gg,\ge_n} \qquad\text{in }S.
 \ee
 This proves that $w_{\gg,\ge_n}$ is a viscosity solution of 
 \bel{Z2'}\BA {ll}
-\myfrac{1}{2}\nabla'\abs{\nabla'w}^2.\nabla'w+\gg\abs{\nabla'w}^4+(2\gg+1)\abs{\nabla'w}^2
+\ge_n w=0\quad\text{in } S\\
\phantom{---------\nabla'W,,,---_{\gg}--,,'}
\displaystyle\lim_{\gr(\gs)\to 0}w(\gs)=\infty,
\EA\ee
and from $(\ref{Z1'})$, 
 \bel{Z2-1}\BA {lll}
 (i)\qquad\qquad &w_{\gg,\ge_n}\geq w_{\gg,\ge_m}\\[2mm]
 (ii)\qquad\qquad &w_{\gg,\ge_n}-w_{\gg,\ge_n}(\gs_0)\geq w_{\gg,\ge_m}-w_{\gg,\ge_m}(\gs_0)
 \\[2mm]
(iii)\qquad\qquad&\ge_nw_{\gg,\ge_n}(\gs_0)\leq \ge_m w_{\gg,\ge_m}(\gs_0).
\EA\ee
Because $\tilde w_{\gg,\ge_n}=w_{\gg,\ge_n}-w_{\gg,\ge_n}(\gs_0)$ is increasing with respect to $n$, locally compact in the topology of local uniform convergence and satisfies
 \bel{Z3}\BA {ll}
-\myfrac{1}{2}\nabla'\abs{\nabla'\tilde w}^2.\nabla'\tilde w+\gg\abs{\nabla'\tilde w}^4+(2\gg+1)\abs{\nabla'\tilde w}^2
+\ge_n\tilde w+\ge_n w'_{\gg,\ge_n}(\gs_0)=0\quad\text{in } S\\
\phantom{---------\nabla'\tilde w,,,---_{k\gg}---------}
\displaystyle\lim_{\gr(\gs)\to 0}\tilde w(\gs)=\infty,
\EA\ee
and since $\ge_n w_{\gg,\ge_n}(\gs_0)\to \gl(\gg,S)$ as $n\to\infty$, we infer that $\tilde w_{\gg}=\lim_{n\to\infty}\tilde w_{\gg,\ge_n}$ is a locally Lipschitz continuous viscosity solution of 
 \bel{Z4}\BA {ll}
-\myfrac{1}{2}\nabla'\abs{\nabla'\tilde w}^2.\nabla'\tilde w+\gg\abs{\nabla'\tilde w}^4+(2\gg+1)\abs{\nabla'\tilde w}^2
+\gl(\gg,S)=0\quad\text{in } S\\
\phantom{----------\nabla'\tilde w,,-------;'}
\displaystyle\lim_{\gr(\gs)\to 0}\tilde w(\gs)=\infty.
\EA\ee

%%%%%%%%%%%%%%%%%%%%%%%%%%%%%%%%%%%%%%%%%%%%%%%%%%%%%%%%%%%%%%%%%%%%%%%%%%%%%%%%%%%%%%%%%%%%%%%%%%%%%%%%%%%%%%%%%%%%%%CASE%2%%%%%%%%%%%%%%%%%%%%%%%%%%%%%%%%%%%%%%%%%%%%%%%%%%%%%%%%%%%%%%%%%%%%%%%%%%%%%%%%%%%%%%%%%%%%%%%%%%%%%%%%%%%%%%%%%%%%%%%%%%%%%%%%%%%%%%%%%%%%%%

\nind{\it Case 2}. If the condition of Step 1 does not hold, for any $\gs_0\in S$ there exist two decreasing sequences $\{\ge_n\}$, $\{\gd_\ell\}$ converging to $0$ and an increasing sequence $\{k_j\}$ tending to infinity, all depending on $\gs_0$, such that 
 \bel{Z7}\BA {ll}
\gl(\gg,S)> \ge_nw_{k_j,\gg,\gd_\ell,\ge_n}(\gs_0)> \ge_mw_{k_m,\gg,\gd_m,\ge_m}(\gs_0)\qquad\forall\,n>m,
\EA\ee
where 
$$\displaystyle\gl(\gg,S)=\lim_{\tiny\BA {lll}n\to\infty\\j\to\infty\\\ell\to \infty\EA}\ge_nw_{k_j,\gg,\gd_\ell,\ge_n}(\gs_0).$$ 
We fix some $\gs_0\in S$. Then $\tilde w_{k_j,\gg,\gd_\ell,\ge_n}=w_{k_j,\gg,\gd_\ell,\ge_n}-w_{k_j,\gg,\gd_\ell,\ge_n}(\gs_0)$ satisfies
 \bel{Z8-0}\BA {ll}
-\gd_\ell\Gd \tilde w-\myfrac{1}{2}\nabla'\abs{\nabla'\tilde w}^2+\gg\abs{\nabla'\tilde w}^4+(2\gg+1)\abs{\nabla'\tilde w}^2
+\ge_n \tilde w+\gl(\gg,S)\geq 0\quad\text{in }\; S_{k_j}\\
\phantom{-----------------,,_{k_j\gg_n\gd_\ell,\ge_n}-_j-'}
\displaystyle\lim_{\gr_{k_j}(\gs)\to 0}\tilde w(\gs)=\infty.
\EA\ee
We introduce the problem 
 \bel{Z7-1}\BA {ll}
-\gd_\ell\Gd Z-\myfrac{1}{2}\nabla'\abs{\nabla'Z}^2+\gg\abs{\nabla'Z}^4+(2\gg+1)\abs{\nabla'Z}^2
+\ge_n Z+\gl(\gg,S)=0\quad\text{in }\; S_{k_j}\\
\phantom{---------------------_\gl'(S)_n}
\displaystyle\lim_{\gr_{k_j}(\gs)\to 0}Z(\gs)=\infty.
\EA\ee
Since $(\ref{Z7-1})$ can be re-written as
 \bel{Z7-2}\BA {ll}-\gd_\ell\Gd Z'-\myfrac{1}{2}\nabla'\abs{\nabla'Z'}^2+\gg\abs{\nabla'Z'}^4+(2\gg+1)\abs{\nabla'Z'}^2
+\ge_n Z'=0\quad\text{in }\; S_{k_j}\\
\phantom{--------------------,'}
\displaystyle\lim_{\gr_{k_j}(\gs)\to 0}Z'(\gs)=\infty,
\EA\ee
with $Z'=Z+\ge_n^{-1}\gl(\gg,S)$, existence is ensured by the approximation by finite boundary data as above. We denote by $Z_{k_j,\gg,\gd_\ell,\ge_n}$ and $Z'_{k_j,\gg,\gd_\ell,\ge_n}$, which coincides actually with $\tilde w_{k_j,\gg,\gd_\ell,\ge_n}$, the solutions of $(\ref{Z7-1})$ and $(\ref{Z7-2})$ obtained by such approximation. Using \rcor{punct-sph} as in Case 1and comparison, we obtain the following estimate
 \bel{Z7-3}\BA {ll}
v_{a,\gg,\ge_n}-\myfrac{\gl(\gg,S)}{\ge_n}\leq Z'_{k_j,\gg,\gd_\ell,\ge_n}-\myfrac{\gl(\gg,S)}{\ge_n}= Z_{k_j,\gg,\gd_\ell,\ge_n}\leq \tilde w_{k_j,\gg,\gd_\ell,\ge_n}\quad\text{in }S_{k_j},
 \EA\ee
where, again $a\in \prt S$ and $v_{a,\gg,\ge_n}$ is the solution of $(\ref{X8})$ with $\ge=\ge_n$ which exists by \rcor{punct-sph}. Now the sequences $\{Z_{k_j,\gg,\gd_\ell,\ge_n}\}_{\ge_n,k_j}$ and  $\{w_{k_j,\gg,\gd_\ell,\ge_n}\}_{k_j}$ are increasing. Letting successively $\gd_\ell\to 0$ and $k_j\to\infty$ we infer that, up to a subsequence, $Z_{k_j,\gg,\gd_\ell,\ge_n}$ converges locally uniformly to some $Z_{\gg,\ge_n}$ and $\tilde w_{k_j,\gg,\gd_\ell,\ge_n}$ converges  locally uniformly to some $\tilde w_{\gg,\ge_n}=w_{\gg,\ge_n}-w_{\gg,\ge_n}(\gs_0)$ which are respectively viscosity solutions of 
  \bel{Z7-4}\BA {ll}
-\myfrac{1}{2}\nabla'\abs{\nabla'Z}^2+\gg\abs{\nabla'Z}^4+(2\gg+1)\abs{\nabla'Z}^2+\gl(\gg,S)
+\ge_n Z=0\quad\text{in } S\\
\phantom{-------------'----_\gl(\gg,S)_n}
\displaystyle\lim_{\gr(\gs)\to 0}Z(\gs)=\infty,
\EA\ee
 and
  \bel{Z7-5}\BA {ll}
-\myfrac{1}{2}\nabla'\abs{\nabla'\tilde w}^2+\gg\abs{\nabla'\tilde w}^4+(2\gg+1)\abs{\nabla'\tilde w}^2
+\ge_nw_{\gg,\ge_n}(\gs_0)+\ge_n \tilde w= 0\quad\text{in } S\\
\phantom{------------------\ge_nw'_{\gg\ge}(\gs_0)}
\displaystyle\lim_{\gr(\gs)\to 0}\tilde w(\gs)=\infty.
\EA\ee
 Furthermore $Z_{\gg,\ge_n}$ and $\tilde w_{\gg,\ge_n}$ are locally bounded in $S$, relatively compact for the local uniform topology and they satisfy
  \bel{Z7-6-0}\BA {ll}
Z_{\gg,\ge_n}\leq \tilde w_{\gg,\ge_n}\quad\text{in }S.
 \EA\ee
 At end, the sequence $\{Z_{\gg,\ge_n}\}_{\ge_n}$ is nondecreasing. Hence, up to a subsequence, $\{\tilde w_{\gg,\ge_n}\}$ converges locally uniformly in $S$ to some $\tilde w_\gg$ which satisfies
   \bel{Z7-7}\BA {ll}\displaystyle
\lim_{\ge_n\to 0}Z_{\gg,\ge_n}=Z_{\gg}\leq \tilde w_{\gg}\quad\text{in }S.
 \EA\ee
 Since 
 $$\lim_{\gr(\gs)\to 0}Z_{\gg,\ge_n}(\gs)=\infty\leq \lim_{\gr(\gs)\to 0}Z_{\gg}(\gs),$$ 
 it follows that $\tilde w_\gg$ is a locally Lipschitz continuous viscosity solution of $(\ref{Z4})$.\smallskip
 
 %%%%%%%%%%%%%%%%%%%%%%%%%%%%%%%%%%%%%%%%%%%%%%%%%%%%%%%%%%%%%%%%%%%%%%%%%%%%%%%%%%%%%%%%%%%%%%%%%%%%%%%%%%%%%%%%%%%%%%%%%%%%%%%%%%%%%%%%%%%%%%%%%%%%%%%%%%%%%%%%%%%%%%%%%%%%%%%%%%%%%%%%%%%%%%%%%%%%%%%%%%%%%%%%%%%%%%%%%%%%%%%%%%%%%%%%%%%%%%%%%%%%%%%%%%%%%%%%%%%%%%%%%%%%%%%%%%%%%%%%%%%%%%%%%%%%%%%%%%%%%%%%%%%%%%%%%%%%%%%%%%%%%%%%%%%%%%%%%%%%%%%%%%%%%%%%%%%%%%%%%%%%%%%%%%%%%%%%%%%%%%%%%%%%%%%%%%%%%%%%%%%%%%%%
\medskip

\nind{\it Step 3: End of the proof.} 
As in the proof of \rprop{existence}, $\gg\mapsto \gl(\gg,S)-\gg$ is a non increasing function of $\gg$. We recall that $\gl(\gg,S^*_\ga)=\gl(\gg,S^*_\ga(a))$. By formula $(\ref{X6})$, for any $\ga>0$ $\lim_{\gg\to 0}\gl(\gg,S^*_\ga)=\infty$. Since
\bel{Y19}\BA {ll}\displaystyle
\gl(\gg,S)-\gg\geq \gl(\gg,S^*_\gp,)-\gg,
\EA\ee
it follows that $\gl(\gg,S^*_\ga)-\gg$ converges to infinity when $\gg$ converges to $0$. Let 
$\ga>0$ such  that $\overline {S_\ga(a)}\subset S$ for some $a\in S$. If 
$$\gg=\myfrac{\gp^2}{4\ga(\gp+\ga)},
$$ 
then $\gl(\frac{\gp^2}{4\ga(\gp+\ga)},S^*_\ga))-\frac{\gp^2}{4\ga(\gp+\ga)}=1$. Since
\bel{Y20}\BA {ll}\displaystyle
\gl\left(\frac{\gp^2}{4\ga(\gp+\ga)},S\right)-\frac{\gp^2}{4\ga(\gp+\ga)}< \gl\left(\frac{\gp^2}{4\ga(\gp+\ga)},S^*_\ga(a)\right)-\frac{\gp^2}{4\ga(\gp+\ga)}=1,
\EA\ee
it follows that
\bel{Y20-1}\BA {ll}\displaystyle
\inf\left\{\gl(\gg,S)-\gg:\gg>0\right\} <1.
\EA\ee
We set
 \bel{Y21}\BA {ll}\displaystyle
\overline \gb_s=\inf\left\{\gg:\gl(\gg,S)-\gg<1\right\}.
\EA\ee
 Let $\{\gg_\gn\}$ be a sequence decreasing to $\overline \gb_s$ when $\gn\to\infty$ and such that $\lim_{\gn\to\infty}\gl(\gg_\gn,S)=\overline \gb_s+1$. As in Step 1 we denote by $w_{k_j,\gg_\gn,\gd_\ell,\ge_n}$  the solution of $(\ref{Y5})$ with $(k_j,\gg_\gn,\gd_\ell,\ge_n)=(k,\gg,\gd,\ge)$. There always holds
 \bel{Y21-0}w_{k_j,\gg_\gm,\gd_\ell,\ge_n}\leq w_{k_j,\gg_\gn,\gd_\ell,\ge_n}\quad\text{if } \gg_\gm\geq \gg_\gn.
\ee
 Again we distinguish two cases \smallskip
 %%%%%%%%%%%%%%%%%%%%%%%%%%%%%%%%%%%%%%%%%%%%%%%%%%%%%%%%%%%%%%%%%%%%%%%%%%%%%%%%%%%%%%%%%%%%%%%%%%%%%%%%%%%%%%%%%%%%%%%%%%%%%%%%%%%%%%%%%%%%%%%%%%%%%%%%%%%%%%%%%%%%%%%%%%%%%%%%
 
 \nind{\it Case 1}. We assume that there exist $\gs_0\in S$ and monotone sequences $\{k_j\}$, $\{\gd_\ell\}$ and $\{\ge_n\}$ such that 
 \bel{Y21-1}\BA {ll}
\ge_nw_{k_j,\gg_\gn,\gd_\ell,\ge_n}(\gs_0)\leq \ge_m w_{k_j,\gg_\gm,\gd_\ell,\ge_m}(\gs_0)\quad\text{for all }\,
m<n\,,\,\gm\leq \gn\,,\;j,\ell\in\BBN.
\EA\ee
 As in Step 2-Case 1 it implies 
  \bel{Y21-2}\BA {ll}
w_{k_j,\gg_\gn,\gd_\ell,\ge_n}-w_{k_j,\gg_\gn,\gd_\ell,\ge_n}(\gs_0)\geq w_{k_j,\gg_\gm,\gd_\ell,\ge_m}- w_{k_j,\gg_\gm,\gd_\ell,\ge_m}(\gs_0)\quad\text{for all }\,
m<n\,,\,\gm<\gn\,,\;j,\ell\in\BBN.
\EA\ee
Letting $\gd_\ell\to 0$ and $k_j\to\infty$ we obtain that the limit function $w_{\gg_\gn,\ge_n}$ satisfies
 \bel{Y22}\BA {ll}
-\myfrac{1}{2}\nabla'\abs{\nabla'w}^2.\nabla'w+\gg_\gn\abs{\nabla'w}^4+(2\gg_\gn+1)\abs{\nabla'w}^2
+\ge_n w=0\quad\text{in } S\\
\phantom{,---------\nabla'W,,,---_{\gg_\gn}--,,'}
\displaystyle\lim_{\gr(\gs)\to 0}w(\gs)=\infty,
\EA\ee
and $\tilde w_{\gg_\gn,\ge_n}=w_{\gg_\gn,\ge_n}-w_{\gg_\gn,\ge_n}(\gs_0)$ is increasing both with respect to $n$ and $\gn$. If $\ge_n\to 0$ we derive that $\tilde w_{\gg_\gn}=\lim_{n\to\infty} \tilde w_{\gg_\gn,\ge_n}$ satisfies $\tilde w_{\gg_\gn}(\gs_0)=0$ and
 \bel{Y23}\BA {ll}-\myfrac{1}{2}\nabla'\abs{\nabla'\tilde w}^2.\nabla'\tilde w+\gg_\gn\abs{\nabla'\tilde w}^4+(2\gg_\gn+1)\abs{\nabla'\tilde w}^2+\gl'(\gg_\gn,S) =0\quad\text{in } S\\
\phantom{,-,\myfrac{1}{2}\nabla'\abs{\nabla'\tilde w'}^2+\gg\abs{}^4+(2\gg+1)\abs{\nabla'}^2-,,---}
\!\displaystyle\lim_{\gr(\gs)\to 0}\tilde w(\gs)=\infty.
\EA\ee
By gradient estimates and since $\tilde w_{\gg_\gn}(\gs_0)=0$, the set of functions $\{\tilde w_{\gg_\gn}\}_\gn$ is relatively compact for the local uniform convergence in $S$. Furthermore $\tilde w_{\gg_\gn}$ is increasing with respect to $\gn$, with limit $\tilde w$. Using $(\ref{Y21})$ and the definition of $\{\gg_\gn\}$, we conclude that 
 \bel{Y23-1}\BA {ll}-\myfrac{1}{2}\nabla'\abs{\nabla'\tilde w}^2.\nabla'\tilde w+\overline\gb_s\abs{\nabla'\tilde w}^4+(2\overline\gb_s+1)\abs{\nabla'\tilde w}^2+\overline\gb_s+1 =0\quad\text{in } S\\
\phantom{,,\nabla\abs{\nabla'\tilde w}^2+\gg\abs{\nabla\tilde w'}^4+(2\gg+1)\abs{\nabla'}^2,',^{,},,--}
\!\displaystyle\lim_{\gr(\gs)\to 0}\tilde w(\gs)=\infty,
\EA\ee
holds in the viscosity sense.\smallskip

\nind{\it Case 2}. We assume that for any $\gs_0\in S$ and $\gn$ there exist two decreasing sequences $\{\ge_n\}$, $\{\gd_\ell\}$ converging to $0$ and an increasing sequence $\{k_j\}$ tending to infinity such that 
 \bel{Y24}\BA {ll}
\overline\gb_s+1> \ge_nw_{k_j,\gg_\gn,\gd_\ell,\ge_n}(\gs_0)> \ge_mw_{k_m,\gg_\gn,\gd_m,\ge_m}(\gs_0)\quad\forall\,n>m,
\EA\ee
where 
$$\overline\gb_s+1=\lim_{\gn\to\infty}\gl(\gg_\gn,S)=\lim_{n\to\infty,j\to\infty,\ell\to \infty}\ge_nw_{k_j,\gg_\gn,\gd_\ell,\ge_n}(\gs_0).$$ 
We follow the ideas in Step 2-Case 2 and consider the problem 
 \bel{Y25}\BA {ll}
-\gd_\ell\Gd Z-\myfrac{1}{2}\nabla'\abs{\nabla'Z}^2+\gg_\gn\abs{\nabla'Z}^4+(2\gg_\gn+1)\abs{\nabla'Z}^2
+\ge_n Z+\overline\gb_s+1=0\quad\text{in } S'_{k_j}\\
\phantom{--------------------,,,-\gb_+^m+}
\displaystyle\lim_{\gr'_{k_j}(\gs)\to 0}Z(\gs)=\infty.
\EA\ee
Since $(\ref{Y25})$ can be re-written as $(\ref{Z7-2})$ with $\gg$ replaced by $\gg_\gn$ and setting 
 $Z'=Z+\ge_n^{-1}(\overline\gb_s+1)$, we have existence and uniqueness of the solution $Z^*_{k_j,\gg_\gn,\gd_\ell,\ge_n}$ (we do not use the previous notation $Z_{k_j,\gg_\gn,\gd_\ell,\ge_n}$ since the constant term is not of the form $\gl(\gg_\gn,S)$). Then ${Z'}^*_{k_j,\gg_\gn,\gd_\ell,\ge_n}=Z^*_{k_j,\gg_\gn,\gd_\ell,\ge_n}+\ge_n^{-1}(\overline\gb_s+1)$ satisfies $(\ref{Z7-2})$ with $\gg$ replaced by $\gg_\gn$.
Then $(\ref{Z7-3})$ is replaced by 
 \bel{Y26}\BA {ll}
v_{a,\gg_\gn,\ge_n}-\myfrac{\overline\gb_s+1}{\ge_n}\leq {Z'}^*_{k_j,\gg_\gn,\gd_\ell,\ge_n}-\myfrac{\overline\gb_s+1}{\ge_n}= Z^*_{k_j,\gg_\gn,\gd_\ell,\ge_n}\leq \tilde w'_{k_j,\gg_\gn,\gd_\ell,\ge_n}\quad\text{in }S_{k_j},
 \EA\ee
where $v_{a,\gg_\gn,\ge_n}$ is as above with obvious modifications. We denote by $\tilde w_{\gg_\gn,\ge_n}=w_{\gg_\gn,\ge_n}-w_{\gg_\gn,\ge_n}(\gs_0)$ the limit, when $\gd_\ell\to 0$ and $k_j\to\infty$, of $\tilde w_{k_j,\gg_\gn,\gd_\ell,\ge_n}=w_{k_j,\gg_\gn,\gd_\ell,\ge_n}-w_{k_j,\gg_\gn,\gd_\ell,\ge_n}(\gs_0)$ and by $Z^*_{\gg_\gn,\ge_n}$ the one of $Z^*_{k_j,\gg_\gn,\gd_\ell,\ge_n}$ under the same conditions. 
They are respective viscosity solutions of
  \bel{Y27}\BA {ll}
-\myfrac{1}{2}\nabla'\abs{\nabla'\tilde w}^2+\gg_\gn\abs{\nabla'\tilde w}^4+(2\gg_\gn+1)\abs{\nabla'\tilde w}^2
+\ge_nw'_{\gg_\gn,\ge_n}(\gs_0)+\ge_n \tilde w= 0\quad\text{in } S\\
\phantom{-------------------\ge_nw'_{\gg_\gn\ge}(\gs_0)}
\displaystyle\lim_{\gr(\gs)\to 0}\tilde w(\gs)=\infty.
\EA\ee
and
\bel{Y28}\BA {ll}
-\myfrac{1}{2}\nabla'\abs{\nabla'Z}^2+\gg_\gn\abs{\nabla'Z}^4+(2\gg_\gn+1)\abs{\nabla'Z}^2+\overline\gb_s+1
+\ge_n Z=0\quad\text{in } S\\
\phantom{,---------------,,-_\gl'(\gg_\gn,S)_n}
\displaystyle\lim_{\gr(\gs)\to 0}Z(\gs)=\infty.
\EA\ee
 Furthermore $Z^*_{\gg_\gn,\ge_n}$ and $\tilde w_{\gg_\gn,\ge_n}$ are locally bounded in $S$, relatively compact for the local uniform topology and they satisfy
  \bel{Z7-6}\BA {ll}
Z^*_{\gg_\gn,\ge_n}\leq \tilde w'_{\gg_\gn,\ge_n}\quad\text{in }S.
 \EA\ee
The sequence $\{Z^*_{\gg_\gn,\ge_n}\}$ is nondecreasing both with respect to $n$ and $\gn$. Therefore the boundary condition is kept. Letting $\ge_n\to 0$ and $\gg_\gn\to \overline\gb_s$ we conclude as in Step 1 that, up to a subsequence 
$\{\gn_s\}$ there exists a locally Lipschitz continuous function $\tilde w$ such that $\tilde w_{\gg_{\gn_s},\ge_n}\to \tilde w$
when $\gn\to \infty$ and $\gn_s \to \overline\gb_s$ succesively, and $\tilde w$ is a viscosity solution of $(\ref{Y23-1})$.\smallskip

%%%%%%%%%%%%%%%%%%%%%%%%%%%%%%%%%%%%%%%%%%%%%%%%%%%%%%%%%%%%%%%%%%%%%%%%%%%%%%%%%%%%%%%%%%%%%%%%%%%%%%%%%%%%%%%%%%%%%%%%%%%%%%%%%%%%%%%%%%%%%%%%%%%%%%%%%%%%%%%%%%%%%%%%%%%%%%%%%%%%%%%%%%%%%%%%%%%%%%%%%%%%%%%%%%%%%%%%%%%%%%%%%%%%%%%%%%%%%%%%%%%%%%%END%%%%%%%%%%%%%%%%%%%%%%%%%%%%%%%%%%%%%%%%%%%%%%%%%%%%%%%%%%%%%%%%%%%%%%%%%%%%%%%%%%%%%%%%%%%%%%%%%%%%%%%%%%%%%%%%%%%%%%%%%%%%%%%%%%%%%%%%%%%%%%%%%%%%%%%%%%%%%%%%%%%

\nind We end the proof by setting $\overline\psi_s=e^{-\overline\gb_s \tilde w}$.\qeda\medskip
%%%%%%%%%%%%%%%%%%%%%%%%%%%%%%%%%%%%%%%%%%%%%%%%%%%%%%%%%%%%%%%%%%%%%%%%%%%%%%%%%%%%%%%%%%%%%%%%%%%%%%%%%%%%%%%%%%%%%%%%%%%%%%%%%%%%%%%%%%%%%%%%%%%%%%%%%%%%%%%%%%%%%%%%%%%%%%%%%%%%%%%%%%%%%%%%%%%%%%%%%%%%%%%%%%%%%%%%%%%%%%%%%%%%%%%%%%%%%%%%%%%%%%%END%%%%%%%%%%%%%%%%%%%%%%%%%%%%%%%%%%%%%%%%%%%%%%%%%%%%%%%%%%%%%%%%%%%%%%%%%%%%%%%%%%%%%%%%%%%%%%%%%%%%%%%%%%%%%%%%%%%%%%%%%%%NEXT THEOREM%%%%%%%%%%%%%%%%%%%%%%%%%%%%%%%%%%%%%%%%%%%%%%%%%%%%%%%%%%%%%%%%%%%%%%%%%%%%%%%%%%%%%%%%%%%%%%%%%%%%%%%%%%%%%%%%%%%%%%%%%%%%%%%%%%%%%%%%%%%%%%%%%%%%%%%%%%%%%%%%%%%%%%%%%%%%%%%%%%%%%%%%%%%%%%%%%%%%%%%%%%%%%%%%%%%%%%%%%%%%%%%%%%%%%%%%%%%%%%%%%%%%%%%%%%%%%%%%%%%%%%%%%%%%%%%%%%%%%%END%%%%%%%%%%%%%%%%%%%%%%%%%%%%%%%%%%%%%%%%%%%%%%%%%%%%%%%%%%%%%%%%%%%%%%%%%%%%%%%%%%%%%%%%%%%%%%%%%%%%%%%%%%%%%%%%%%%%%%%%%%%%%%%%%%%%%%%%%%%%%%%%%%%%%%%%%%%%%%%%%%%

{\it Mutatis mutandis} in the above proof, one can obtain an existence result of a separable positive regular infinity harmonic function in $C_S$ vanishing on $\prt C_S$.

\bth{D'} Assume $S\subsetneq S^{N-1}$ is any domain. Then there exist $\overline\gm_s>0$ and a positive function $\overline \gw_s$ in $C(\overline S)$, locally Lipschitz continuous in $S$ and vanishing on $\prt S$, such that the function 
\bel{I9-2}\BA {ll}\displaystyle
\overline u_{s,{\scriptscriptstyle-}}(r,\gs)=r^{\overline\gm_s}\overline \gw_s(\gs),
\EA\ee
is infinity harmonic in $C_{S}$ and vanishes on $\prt C_{S}$.
\es

\subsection{Proof of Theorems E and F}

\nind{\it Proof of Theorem E}. {\it Step 1: Existence of $\underline\gb_s$}. The proof follows the one of Theorem D, hence we indicate only the main streamlines. We consider a decreasing sequence of smooth spherical domains $\{S'_k\}$ such that 
$$\displaystyle S\subset S'_{k+1}\subset\overline S'_{k+1}\subset S'_{k}\,\text{ and }\;\rm{int}\left(\bigcap_k S'_k\right)=S.
$$
Such a sequence of domains $\{S'_k\}$ exists since $\prt S=\prt\overline S^c$. To each domain we associate  the positive exponent $\gb'_k:=\gb_{s'_k}$ and the corresponding spherical p-harmonic function 
$\psi'_k:=\psi_{s'_k}$, defined in $S'_k$, such that $\psi'_k(\gs_*)=1$ for some $\gs_*\in S_1$, so that the function 
 $u'_{k}(r,.)=r^{-\gb'_{k}}\psi'_{k}$ is $p$-harmonic in 
 $C_{S'_k}$ and vanish on $\prt C_{S'_k}\setminus\{0\}$. For $\gg,\gd,\ge>0$, we denote by $w'_{k,\gg,\gd,\ge}$  the solution of
\bel{Y6}\BA {ll}
-\gd\Gd w'-\myfrac{1}{2}\nabla'\abs{\nabla'w'}^2.\nabla'w'+\gg\abs{\nabla'w'}^4+(2\gg+1)\abs{\nabla'w'}^2+\ge w'=0\quad\text{in } S'_k\\[2mm]\phantom{.\gd\Gd w'-\myfrac{1}{2}\nabla'\abs{\nabla'w'}^2+\gg\abs{\nabla'w'}^4+(2\gg+1)\abs{\nabla'w}^2}
\displaystyle\lim_{\gr'_k(\gs)\to 0}w'(\gs)=\infty,
\EA\ee
where $\gr'_k(.)=\dist(.,\prt S'_k)$. By the maximum principle all the functions 
$w'_{\ell,\gg,\gd,\ge}$ is positive and  the following comparison relations hold. Estimates $(\ref{Y6+})$ are valid the main difference being the fact that $w'_{\ell,\gg,\gd,\ge}\leq w_{k,\gg,\gd,\ge}$ in $S_k$ for $k,\ell>0$ and that  the mapping $k\mapsto w'_{k,\gg,\gd,\ge}$ is increasing. Similarly  $(\ge,\gg)\mapsto w'_{\ell,\gg,\gd,\ge}$ is nonincreasing. The gradient estimate $(\ref{Y7})$ holds for $w'_{k,\gg,\gd,\ge}$, provided 
$\gr_k(\gs)$ be replaced by $\gr'_k(\gs)=\dist (\gs,\prt S'_k)$.
Moreover, similarly to in $(\ref{G1-0})$,
\bel{Y8+'}\BA {llll}
\qquad\; -\myfrac{1}{\gg}\ln\gr'_k(\gs)-\myfrac{M'_k}{\ge}\leq  w'_{k,\gg,\gd,\ge}(\gs)\leq -\myfrac{1}{\gg}\ln\gr'_k(\gs)+\myfrac{M'_k}{\ge}\qquad\;\forall \gs\in S'_k.
\EA\ee
When $\gd\to 0$, $w'_{k,\gg,\gd_n,\ge}\to w'_{k,\gg,\ge}$ locally uniformly in $S'_k$ to some function $w'_{k,\gg,\ge}$ which satisfies $(\ref{Y6+})$, the modified gradient estimate $(\ref{Y7})$ (expressed with $\gr_k$ replaced by $\gr'_k$) and $(\ref{Y8+'})$ and is a viscosity solution of 
\bel{Y10+'}\BA {ll}
-\myfrac{1}{2}\nabla'\abs{\nabla'w'}^2.\nabla'w'+\gg\abs{\nabla'w'}^4+(2\gg+1)\abs{\nabla'w'}^2+\ge w'=0\quad\text{in } S'_k\\[2mm]\phantom{.\myfrac{1}{2}\nabla'\abs{\nabla'w'}^2+\gg\abs{\nabla'w'}^4+(2\gg+1)\abs{\nabla'w}^2}
\displaystyle\lim_{\gr'_k(\gs)\to 0}w'(\gs)=\infty.
\EA\ee
When $k\to\infty$ $w'_{k,\gg,\ge}\uparrow w'_{\gg,\ge}$ which is a nonincreasing function of $\ge$ and $\gg$.\smallskip

The proof of the boundary blow-up introduces two cases: either there exists $\gs_0\in S$, two decreasing sequences $\{\ge_n\}$, $\{\gd_\ell\}$ converging to $0$ and an increasing sequence $\{k_j\}$ tending to infinity with the property that
 \bel{Z1+'}\BA {ll}
\ge_nw'_{k_j,\gg,\gd_\ell,\ge_n}(\gs_0)\leq \ge_m w'_{k_j,\gg,\gd_\ell,\ge_m}(\gs_0)\quad\text{for all }\,
m<n\,,\;j,\ell\in\BBN.
\EA\ee
Or for any $\gs_0\in S$ there exist two decreasing sequences $\{\ge_n\}$, $\{\gd_\ell\}$ converging to $0$ and an increasing sequence $\{k_j\}$ tending to infinity, all depending on $\gs_1$, such that 
 \bel{Z7+'}\BA {ll}
\gl'(\gg,S)> \ge_nw'_{k_j,\gg,\gd_\ell,\ge_n}(\gs_0)> \ge_mw'_{k_m,\gg,\gd_m,\ge_m}(\gs_0)\qquad\forall\,n>m,
\EA\ee
where 
$$\displaystyle\gl'(\gg,S)=\lim_{\tiny\BA {lll} n\to\infty,\\j\to\infty,\\\ell\to \infty\EA}\ge_nw'_{k_j,\gg,\gd_\ell,\ge_n}(\gs_1).$$ 
In the first case for any $a\in\prt S$, there exists a sequence $\{a_j\}$ such that $a_j\in S'_{k_j}$ converging to $a$ (such a sequence exits since $\prt S=\prt \bar S^c$). Then 
 \bel{Z1''+'}
v_{a_j,\gg,\ge_n} \leq w'_{k_j,\gg,\ge_n} \qquad\text{in }S_{k_j}
 \ee
Since $v_{a,\gg,\ge_n}$ is obtained from $v_{a_j,\gg,\ge_n}$ by an orthogonal transformation on $S^{N-1}$, we derive
 \bel{Z1'''}
v_{a,\gg,\ge_n} \leq w'_{\gg,\ge_n} \qquad\text{in }S.
 \ee
 This proves that $w'_{\gg,\ge_n}$ is a viscosity solution of 
 \bel{Z2'}\BA {ll}
-\myfrac{1}{2}\nabla'\abs{\nabla'w}^2.\nabla'w+\gg\abs{\nabla'w}^4+(2\gg+1)\abs{\nabla'w}^2
+\ge_n w=0\quad\text{in } S\\
\phantom{---------\nabla'W,,,---_{\gg}--,,'}
\displaystyle\lim_{\gr(\gs)\to 0}w(\gs)=\infty.
\EA\ee
The proof in the second case is the same as in Theorem D, just replacing $v_{a,\gg,\ge_n}$ by $v_{a_j,\gg,\ge_n}$ in $(\ref{Z7-3})$ which becomes
 \bel{Z7-3'}\BA {ll}
v_{a_{k_j},\gg,\ge_n}-\myfrac{\gl(\gg,S)}{\ge_n}\leq Z'_{k_j,\gg,\gd_\ell,\ge_n}-\myfrac{\gl(\gg,S)}{\ge_n}= Z_{k_j,\gg,\gd_\ell,\ge_n}\leq \tilde w_{k_j,\gg,\gd_\ell,\ge_n}\quad\text{in }S_{k_j},
 \EA\ee
where $Z'_{k_j,\gg,\gd_\ell,\ge_n}$ and $Z'_{k_j,\gg,\gd_\ell,\ge_n}$ are defined accordingly. This implies again that the limit $w'_{\gg,\ge_n}$ of 
$\{w'_{k_j,\gd_\ell,\ge_n}\}$ when $j\to\infty$ is a viscosity solution of $(\ref{Z2'})$.\\
The proof of the existence of some $\underline\gb_s>0$ such that $\gl'(\underline\gb_s,S)=1+\underline\gb_s$ follows the same dichotomy.\medskip

\nind {\it Step 2: Comparison of exponents}. Since $w'_{k,\gg,\gd,\ge}\leq w_{k,\gg,\gd,\ge}$ it follows that 
$\ge w'_{k,\gg,\gd,\ge}(\gs_0)\leq \ge w_{k,\gg,\gd,\ge}(\gs_0)$ (it is always possible to choose the same $\gs_0$ in order to defined the ergodic constant), hence $\gl'(\gg,S)\leq \gl(\gg,S)$ and finally $\underline\gb_s\leq \overline\gb_s$ by monotonicity. Assume now that 
$u(r,\gs)=r^{-\gb}\psi(\gs)$ is a positive infinity harmonic function in $C_S$ which vanishes on $\prt C_S\setminus\{0\}$. We proceed by conradiction in assuming that $\gb>\overline \gb_s$. Hence $\gb>\gb_k:=\gb_{s_k}$ for $k$ large enough. We set 
$$\gf=\psi^\gth\quad\text{where }\;\;\gth=\myfrac{\gb_k}{\gb}<1.
$$
Then 
$$\BA{lll}\nabla\phi=\gth\psi^{\gth-1}\nabla\psi\;,\;\,|\nabla\phi|^2=\gth^2\psi^{2(\gth-1)}|\nabla\psi|^2 \;,\;\,
|\nabla\phi|^2\phi=\gth^2\psi^{3(\gth-1)+1}|\nabla\psi|^2,\\ [1mm]
\nabla|\nabla\phi|^2.\nabla\phi=\gth^3\psi^{3(\gth-1)}\nabla|\nabla\psi|^2.\nabla\psi+2\gth^3(\gth-1)\psi^{3(\gth-1)-1}|\nabla\psi|^4
\EA$$
If we denote 
$$\CL_k\phi=-\myfrac{1}{2}\nabla'\abs{\nabla'\phi}^2.\nabla'\phi-\gb_k(2\gb_k+1)\abs{\nabla'\phi}^2\phi-\gb_k^3(\gb_k+1)\phi^3,$$
then
 \bel{XZ1}\BA {ll}
 \CL_k\phi=\gth^3\psi^{3(\gth-1)}\left[\left(\gb(2\gb+1)-\myfrac{\gb_k}{\gth}(2\gb_k+1)\right)\psi|\nabla\psi|^2+(1-\gth)\psi^{-1}|\nabla\psi|^4\right.\\[4mm]\phantom{-----------------}
 \left.  +\left(\gb^3(\gb+1)-\myfrac{\gb^3_k}{\gth^3}(\gb_k+1)\psi^3\right)\right]
 \\[4mm]
 \phantom{ \CL_k\phi}
 =\myfrac{\gth^3\psi^{3(\gth-1)-1}(\gb-\gb_k)}{\gb}\left(|\nabla\psi|^2+\gb^2\psi^2\right)^2>0.
\EA\ee
Hence the function $W=-\frac{1}{\gb_k}\ln\phi$ satisfies in the viscosity sense
$$-\myfrac{1}{2}\nabla'\abs{\nabla'W}^2.\nabla'W+\gb_k|\nabla W|^4+(2\gb_k+1)|\nabla W|^2+\gb_k+1< 0\quad\text{in }S_k$$
and is bounded in $S_k$. By comparison between viscosity solutions, $W$ is smaller than $w_k:=-\frac{1}{\gb_k}\ln\psi_{s_k}$ where $\psi_{s_k}$ is a spherical infinity harmonic function in $S_k$ (see the proof of Theorem D-Step 1). But we can replace $W$ by $W+n$ for any $n\in\BBN_*$. This is a contradiction, hence $\gb\leq \overline\gb_s$. In the same way $\gb\geq \underline\gb_s$, which ends the proof.\qeda\medskip

\nind{\it Proof of Theorem F.} If $\prt S$ is Lipschitz and satifies the interior sphere condition, it is possible to construct a bounded Lipschitz subdomain $\Gth$ of $C_S$ satisfying the interior sphere condition with the following additional properties:
 \bel{XZ2}\BA {llll}
 \Gth\subset C_S\cap \left(B_2\setminus\overline B_{\frac12}\right)\\[2mm]
\Gth_{\frac32,\frac23}:= \Gth\cap \left(B_\frac32\setminus\overline B_{\frac23}\right)
 =C_S\cap \left(B_\frac32\setminus\overline B_{\frac23}\right)=\left\{(r,\gs):\frac{2}{3}<r<\frac32,\gs\in S \right\}.
\EA\ee
and we define the lateral boundary of $\Gth$ by  
$$\prt_\ell\Gth_{\frac32,\frac23}=\prt C_S\cap \left(B_\frac32\setminus\overline B_{\frac23}\right).$$
Assume now that $u(r,\gs)=r^{-\gb}\psi(\gs)$ and $u'(r,\gs)=r^{-\gb'}\psi'(\gs)$, with $\gb,\gb'>0$, are nonnegative, infinity harmonic in $C_S$ and vanish on $\prt C_S\setminus\{0\}$. Since $u$ and $u'$ vanish on $\prt_\ell\Gth_{\frac32,\frac23}$, it follows from \cite[Th 1.1]{Juu2}, that there exists a constant $c_1=c_1(N,\Gth)>0$ such that for any $z\in \prt_\ell\Gth_{\frac32,\frac23}\cap\prt S$ there exists a constant $\gd_z\in (0,\frac{1}{4})$ 
such that 
 \bel{XZ3}\BA {llll}
\myfrac{u(x)}{u'(x)}\leq c_1\myfrac{u(y)}{u'(y)}\qquad\forall (x,y)\in \Gth_{\frac32,\frac23}\cap B_{\gd_z}(z).
\EA\ee
By compactness of $\prt_\ell\Gth_{\frac32,\frac23}\cap\prt S$ the constant $\gd_z$ is actually independent of $z$ and denoted by $\gd^*$. We use now the standard Harnack inequality for infinity harmonic functions in  $\Gth$ (see e.g. \cite{Bat1}) to derive the existence of $c_2=c_2(N,\gd^*)>0$ such that 
 \bel{XZ4}\BA {llll}
\myfrac{1}{c_2}\leq \myfrac{u(x)}{u(y)},\myfrac{u'(x)}{u'(y)}\leq c_2\qquad\forall (x,y)\in \Gth_{\frac32,\frac23}\,\text{ s.t. }\inf\left\{\dist (x,\prt \Gth),\dist (y,\prt \Gth)\right\}\geq \gd^*.
\EA\ee
Taking $x=(1,\gs)$ and $y=(1,\gs')$ in $( \ref{XZ3})$, we derive from $(\ref{XZ3})$, $(\ref{XZ4})$ that
$$\myfrac{\psi(\gs)}{\psi'(\gs)}\leq c_3\myfrac{\psi(\gs')}{\psi'(\gs')}\qquad\forall (\gs,\gs')\in S.
$$
This implies that for a fixed $\gs_0\in S$, one has
 \bel{XZ5}\BA {llll}
\myfrac{1}{c_3}\myfrac{\psi(\gs_0)}{\psi'(\gs_0)}\leq \myfrac{\psi(\gs)}{\psi'(\gs)}\leq c_3\myfrac{\psi(\gs_0)}{\psi'(\gs_0)}\qquad\forall \gs\in S.
\EA\ee
We proceed as is Step 1, assuming $\gb>\gb'$ and defining
$$\phi_*=\psi^\gth\quad\text{where }\;\;\gth_*=\myfrac{\gb'}{\gb}<1
$$
Then
$$-\myfrac{1}{2}\nabla'\abs{\nabla'\phi_*}^2.\nabla'\phi_*-\gb'(2\gb'+1)\abs{\nabla'\phi_*}^2\phi_*-\gb'^3(\gb'+1)\phi_*^3> 0.
$$
The function $W_*=-\frac{1}{\gb'}\ln\phi_*$ satisfies the inequation
 \bel{XZ6}-\myfrac{1}{2}\nabla'\abs{\nabla'W_*}^2.\nabla'W_*+\gb'\nabla W_*|^4+(2\gb'+1)|\nabla W_*|^2+\gb'+1< 0\quad\text{in }S,
 \ee
while $w'=-\frac{1}{\gb'}\ln\psi'$ is a solution of the associated equation. 
From $(\ref{XZ5})$, $\psi'(\gs)=o(\phi(\gs))$ when $\gr(\gs):=\dist (\gs,\prt S)\to 0$. Hence 
 \bel{XZ7}w'-W_*=-\frac{1}{\gb'}\ln\left(\myfrac{\psi'}{\phi_*}\right)\to \infty\quad\text{when }\,\gr(\gs)\to 0.
\ee
By comparison there holds $w'\geq W_*$. Since for any $n\in \BBN_*$, the function $W_{*,n}=W_*+n$ satisfies $(\ref{XZ6})$ and $(\ref{XZ7})$, it follows that $w'\geq W_{*,n}$, contradiction. Hence $\gb\leq \gb'$, which ends the proof.\qeda
%%%%%%%%%%%%%%%%%%%%%%%%%%%%%%%%%%%%%%%%%%%%%%%%%%%%%%%%%%%%%%%%%%%%%%%%%%%%%%%%%%%%%%%%%%%%%%%%%%%%%%%%%%%%%%%%%%%%%%%%%%%%%%%%%%%%%%%%%%%%%%%%%%%%%%%%%%%%%%%%%%%%%%%%%%%%%%%%%%%%%%%%%%%%%%%%%%%%%%%%%%%%%%%%%%%%%%%%%%%%%%%%%%%%%%%%%%%%%%%%%%%%%%%%%%%%%%%%%%%%%%%%%%%%%%%%%%%%%%%%%%%%%%%%%%%%%%%%%%%%

\nind Marie-Fran\c{c}oise Bidaut-V\'eron \\
Laboratoire de Math\'ematiques et Physique Th\'eorique\\
Universit\'e Francois Rabelais, Tours, France\\
{\it veronmf@univ-tours.fr}\\

\nind Marta Garcia Huidobro\\
Departamento de Mathematic\`a\\
Pontificia Universit\`a Catolica, 
Santiago, Chile\\
{\it martaghc@gmail.com}\\

\nind Laurent V\'eron \\
Laboratoire de Math\'ematiques et Physique Th\'eorique\\
Universit\'e Francois Rabelais, Tours, France\\
{\it veronl@univ-tours.fr}
%%%%%%%%%%%%%{\it }{\bf }%%%%%%%%%%%%%%%%%%%%%%
%%%%%%%%%%%%%%%%%%%%%%%%%%%%%%%%%%%
%%%%%%BEGINING OF THE ARTICLE%%%%%%%%%%%%%%
%%%%%%%%INTRODUCTION%%%%%%%%%%%%%%%%%%%%%%%%%%%%%%%%%%%%%%%%%%%%%%%%%%%%%
%%%%%%%%%%%%%%%%%%%%%%%%%%%%%%%%%%%

%Furthermore%Furthermore%Furthermore
%Furthermore%Furthermore%Furthermore
%Furthermore%Furthermore%Furthermore
 %%END DOCUMENT%%%%%%%%%%%%%%%%%%%%%%%%

\begin{thebibliography}{99}

\bibitem{Aro}G. Arronson: Extension of functions satisfying Lipschitz conditions. {\it Ark. Mat.} {\bf 6} (1967), 551-561.

\bibitem{Aro1} G. Arronson: On the equation $u_xu_{xx}+2u_{xy}u_xu_y+u_yu_{yy}=0$. {\it  Ark. Mat.} {\bf 7} (1968), 395-425.

\bibitem{BanMar} C. Bandle and M. Marcus: Dependence of Blowup Rate of Large Solutions
of Semilinear Elliptic Equations, on the Curvature
of the Boundary. {\it Complex Variables},
{\bf 49} (2004), 555-570.

\bibitem{Bat} T. Bhattacharya: A note on non-negative singular infinity-harmonic functions in the half-space.
{\it Rev. Mat. Complut. } {\bf 18} (2005), 377-385.

\bibitem{Bat1} T. Bhattacharya: A boundary Harnack principle for infinity-Laplacian and some related results. {\it Boundary value problems 2007}, 1-17.

\bibitem{EvSm} E. C. Evans and Ch. K. Smart: Everywhere differentiability of infinity harmonic
functions. {\it Calc. Var. Part. Diff. Eq.} {\bf 42} (2011), 289-299.

\bibitem{GkVe} K. Gkikas and L. V\'eron: The spherical p-harmonic eigenvalue problem in non-smooth domains. {\it J. Funct. Anal.}, https://doi.org/10.1016/j.jfa.2017.07.012 (2018), in press.

\bibitem{Kro} I. N. Krol: The behaviour of the solutions of a certain quasilinear equation near zero cusps of
the boundary. {\it Proc. Steklov Inst. Math.} {\bf 125} (1973), 140-146.


\bibitem{Kat} N. Katzourakis: An Introduction to Viscosity Solutions for Fully Nonlinear PDE with Applications to Calculus of Variations in $L^\infty$. Springer Briefs in Mathematics (2015), Springer-Verlag.

\bibitem{KiVe} S. Kichenassamy and L. V\'eron: Singular solutions of the $p$-Laplace equation. {\it Math. Ann.} {\bf 275} (1986), 599-615.

\bibitem{Juu} P. Juutinen:  Principal eigenvalue of a very badly
degenerate operator and applications. {\it J. Diff. Equ.} {\bf 236} (2007), 532-550.

\bibitem{Juu2} P. Juutinen:  The boundary Harnack inequality for infinity harmonic functions in Lipschitz domains satisfying the interior ball condition. {\it Nonlinear Anal.} {\bf 69} (2007), 1941-1944.

\bibitem{LaLi} J. M. Lasry and P. L. Lions: Nonlinear elliptic equations with singular boundary conditions and
stochastic control with state constraints. I. The model problem. {\it Math. Ann.} {\bf 283} (1989), 583-630.

\bibitem{LeNy} J. Lewis, K. Nystr\"{o}m:
Boundary behavior and the Martin boundary problem for $p$-harmonic functions in Lipschitz domains.
{\it Ann. of Math.} {\bf 172}  (2010), 1907-1948.


\bibitem{PoVe} A. Porretta and L. V\'eron: Separable p-harmonic functions in a cone and
related quasilinear equations on manifolds. {\it J. Europ. Math. Sci.} {\bf 11} (2009), 1285-1305.

\bibitem{Tolk} P. Tolksdorf: On the Dirichlet problem for quasilinear equations in domains with conical
boundary points. {\it Comm. Part. Diff. Eq.} {\bf  8} (1983), 773-817.

\bibitem{Tolk 2} P. Tolksdorf: Regularity for a more general class of quasilinear elliptic equations. {\it J. Diff. Eq.} {\bf  51} (1984), 126-150.

\end{thebibliography}
\end {document}